\documentclass[a4paper]{article}

\usepackage{latexsym}
\usepackage{amsmath}
\usepackage{amssymb}
\usepackage{amsthm}
\usepackage{hyperref}
\usepackage{caption}

\newtheorem{Theorem} {Theorem} [section]
\newtheorem{Corollary} [Theorem] {Corollary}
\newtheorem{Proposition} [Theorem] {Proposition}
\newtheorem{Lemma} [Theorem] {Lemma}
\newtheorem{Conjecture} [Theorem] {Conjecture}
\newtheorem*{Conjecture*} {Conjecture}
\newcommand{\Proof}{ \noindent{\bf Proof.}\quad }
\newcommand{\Proofn}{ \noindent{\bf Proof }}
\newcommand{\Remark}{ \noindent{\bf Remark.}\quad }

\newtheoremstyle{named}{}{}{\itshape}{}{\bfseries}{}{ }{{#1}}
\theoremstyle{named}
\newtheorem{nTheorem}{Theorem}

\newtheorem{nCorollary}{Corollary}

\begin{document}

\newcommand{\maysplit}{\discretionary{}{}{}}
\newcommand{\Ff}{{\mathbb F}}
\newcommand{\Zz}{{\mathbb Z}}
\newcommand{\bnom}[2]{\genfrac{[}{]}{0pt}{}{#1}{#2}}
\newcommand{\qnom}[2]{\genfrac{[}{]}{0pt}{}{#1}{#2}}
\newcommand{\txtfrac}[2]{\frac{\textstyle {#1}}{\textstyle {#2}}}



\title{The smallest eigenvalues of Hamming graphs, Johnson graphs
and other distance-regular graphs with classical parameters}


\author{Andries E. Brouwer\footnote{Department of Mathematics, Technische Universiteit Eindhoven, Eindhoven, Netherlands {\tt aeb@cwi.nl}.},
Sebastian M. Cioab\u{a}\footnote{Department of Mathematical Sciences, University of Delaware, Newark, DE 19716-2553, USA {\tt cioaba@udel.edu}.} ,\\
Ferdinand Ihringer\footnote{\hbadness=1500 Einstein Institute of Mathematics, Hebrew University of Jerusalem, Jerusalem, Israel. Department of Pure Mathematics and Computer Algebra, Ghent University, Belgium. {\tt Ferdinand.Ihringer@gmail.com}.}\, and
Matt McGinnis\footnote{Department of Mathematical Sciences, University of Delaware, Newark, DE 19716-2553, USA {\tt mamcginn@udel.edu}.}}

\maketitle

\begin{abstract}
We prove a conjecture by Van Dam \& Sotirov on the smallest eigenvalue
of (distance-$j$) Hamming graphs and a conjecture by Karloff on the
smallest eigenvalue of (distance-$j$) Johnson graphs.
More generally, we study the smallest eigenvalue and the second largest
eigenvalue in absolute value of the graphs of the relations of classical
$P$- and $Q$-polynomial association schemes.
\end{abstract}
%
%
%
%

\section{Introduction}

In this paper we study the smallest eigenvalue as well as the second largest one
in absolute value of the adjacency matrix of several important families of graphs,
all belonging to the classical $P$- and $Q$-polynomial association schemes \cite[Chapter 6]{BCN}.

The most well-known example of a $P$-polynomial association scheme
is the Hamming scheme. We investigate the eigenvalues of the graphs 
that have the vectors in $\Ff_q^d$ as vertices 
and two vertices are adjacent if they have Hamming distance $j$. 
The smallest eigenvalues are important for determining the max-cut of 
certain graphs in the Hamming scheme. These graphs provide examples where the performance ratio of the Goemans-Williamson algorithm is tight \cite{AS}.
The smallest eigenvalues are also used for 
determining the max-$k$-cut \cite{vDS} and the chromatic number of the graphs 
in the Hamming scheme \cite{vDS}.

The second important scheme belonging to the family of $P$-polynomial 
association schemes is the Johnson scheme.
Here the vertices are the $d$-subsets of $\{ 1, 2, \ldots, n \}$.
We investigate the eigenvalues of the graph where two $d$-sets are adjacent
if they differ in exactly $j$ elements. 
As for the Hamming scheme, these graphs provide examples 
for which the performance ratio of the Goemans-Williamson algorithm 
is tight and their smallest eigenvalues are central for determining their max-cuts \cite{Karloff}. 
These graphs are also important for investigating subsets 
with exactly one forbidden intersection, a variation of the classical 
Erd\H{o}s-Ko-Rado theorem due to Frankl and F\"{u}redi \cite{FF}.

The other graphs under investigation are 
Grassmann graphs, dual polar graphs,
and various forms graphs, most prominently the bilinear forms graphs. 
Again, the smallest eigenvalues can be used to investigate the max-cuts
and intersecting families in these graphs.
The $P$-polynomial graphs obtain their importance from various applications. For example, Grassmann graphs 
are of interest due to their applications in network coding theory \cite{Silva2008}
and their role in the recent proof of the $2$-to-$2$-games conjecture \cite{Khot2018}.

In the following we give a short summary of our main results on the specific families.

\subsection{Hamming graphs}\label{subsec:intro_Hamming}

Let $q \ge 2, d \ge 1$ be integers. Let $Q$ be a set of size $q$.
The Hamming scheme $H(d,q)$ is the association scheme with vertex set
$Q^d$, and as relation the Hamming distance. The $d+1$ relation graphs
$H(d,q,j)$, where $0 \le j \le d$, have vertex set $Q^d$, and two
vectors of length $d$ are adjacent when they differ in $j$ places.

The eigenmatrix $P$ of $H(d,q)$ has entries $P_{ij} = K_j(i)$, where
$$
K_j(i) = \sum_{h=0}^j (-1)^h (q-1)^{j-h} \binom{i}{h} \binom{d-i}{j-h} .
$$
The eigenvalues of the graph $H(d,q,j)$ are the numbers in column $j$ of $P$,
so are the numbers $K_j(i)$, $0 \le i \le d$.
The graph $H(d,q,j)$ is regular of degree $K_j(0) = (q-1)^j \binom{d}{j}$,
and this is the largest eigenvalue.
Motivated by problems in semidefinite programming related to the max-cut
of a graph, Van Dam \& Sotirov \cite{vDS} conjectured

\begin{Conjecture}\label{conj:vandam_sotirov}
Let $q \ge 2$ and $j \ge d - \frac{d-1}{q}$ where $j$ is even when
$q = 2$. Then the smallest eigenvalue of $H(d,q,j)$ is $K_j(1)$.
\end{Conjecture}

Alon \& Sudakov \cite{AS} proved this for $q=2$ and $d$ large and $j/d$ fixed.
Dumer \& Kapralova \cite[Cor.~10]{DK}, proved this for $q=2$ and all $d$.
Here we settle the full conjecture.

In most cases $K_j(1)$ is not only the smallest eigenvalue,
but also the second largest eigenvalue in absolute value.
The only exception is the case $d=4$, $q=3$:
the $P$-matrix of $H(4,3)$ is
$$
P = \left(\begin{array}{ccccc}
 1 &   8 &  24 & 32 & 16 \\
 1 &   5 &   6 & -4 & -8 \\
 1 &   2 &  -3 & -4 &  4 \\
 1 &  -1 &  -3 &  5 & -2 \\
 1 &  -4 &   6 & -4 &  1
\end{array}\right)
$$
and the eigenvalues of $H(4,3,3)$ are $-4$, $5$ and $32$.

\medskip

The binary case was already settled by Dumer \& Kapralova.
We give a short and self-contained proof.

\begin{Theorem}\label{binary} {\rm (\cite[Cor.~10]{DK})} Let $q = 2$.

(i) If $j \ne d/2$, then $|K_j(i)| \le |K_j(1)|$
for all $i$, $1 \le i \le d-1$.

(ii) If $j = d/2$, then $K_j(1) = 0$ and $|K_j(i)| \le |K_j(2)|$
for all $i$, $1 \le i \le d-1$.
\end{Theorem}

\begin{Corollary}\label{binarycor}
Let $q=2$ and $j \ge (d+1)/2$.

(i) One has $K_j(1) \le K_j(i)$ for all $i$, $0 \le i \le d-1$.

(ii) One has $K_j(1) \le K_j(d)$ if and only if $j$ is even or $j = d$.
\end{Corollary}

The nonbinary case is settled here.

\begin{Theorem}\label{nonbinary}
Let $q \ge 3$ and $d - \frac{d-1}{q} \le j \le d$.

(i) One has $K_j(1) \le K_j(i)$ for all $i$, $0 \le i \le d$.

(ii) One has $|K_j(i)| \le |K_j(1)|$ for all $i \ge 1$,
unless $(q,d,i,j) = (3,4,3,3)$.
\end{Theorem}

\subsection{Johnson graphs} \label{subsec:intro_Johnson}
The Johnson graphs $J(n,d)$ are the graphs with as vertices the $d$-subsets
of a fixed $n$-set, adjacent when they meet in a $(d-1)$-set.
W.l.o.g. we assume $n \ge 2d$ (since $J(n,d)$ is isomorphic to $J(n,n-d)$),
and then these graphs are distance-regular of diameter $d$.
The eigenmatrix $P$ has entries $P_{ij} = E_j(i)$, where
$$
E_j(i) = \sum_{h=0}^i (-1)^{i-h} \binom{i}{h} \binom{d-h}{j} \binom{n-d-i+h}{n-d-j}.
$$

For $0 \le j \le d$, the distance-$j$ graphs $J(n, d, j)$ of
the Johnson graph $J(n, d)$ are the graphs with the same vertex set
as $J(n, d)$, where two vertices are adjacent when they have
distance $j$ in $J(n, d)$, that is, when they meet in a \mbox{$(d-j)$}-set.
For $j=d$ this graph is known as the Kneser graph $K(n,d)$.
Motivated by problems in semidefinite programming related to the max-cut
of a graph, Karloff \cite{Karloff} conjectured in 1999 the following:

\begin{Conjecture}
Let $n = 2d$ and $j > d/2$.
Then the smallest eigenvalue of $J(n,d,j)$ is $E_j(1)$.
\end{Conjecture}

Here we prove this conjecture (Corollary \ref{Kconj}),
and more generally determine precisely in which cases
$E_j(1)$ is the smallest eigenvalue of $J(n,d,j)$
(Theorem \ref{smallestj1}).

\subsection{Graphs with classical parameters}\label{subsec:intro_DRG}
For general information on distance-regular graphs, see \cite{BCN}.
In \cite[\S6.1]{BCN}, {\em graphs with classical parameters}
$(d,b,\alpha,\beta)$ are defined as distance regular graphs
of diameter $d$ with parameters given by certain expressions
in $d,b,\alpha,\beta$ (see Section \ref{classparams} for details).

The concept of graphs with classical parameters unifies a number
of families of distance-regular graphs, such as the Hamming graphs,
Johnson graphs, Grassmann graphs, dual polar graphs, bilinear forms graphs,
etc.

\medskip\noindent
\begin{tabular}{@{}c@{~~~}c@{~~~}c@{~~~}c@{~~~}l@{}}
$d$ & $b$ & $\alpha$ & $\beta$ & family \\
\hline
$d$ & 1 & 0 & $q-1$ & Hamming graphs $H(d,q)$ \\
$d$ & 1 & 1 & $n-d$ & Johnson graphs $J(n,d)$, $n \ge 2d$ \\
$d$ & $q$ & $q$ & $q\qnom{n-d}{1}$ & Grassmann graphs $G_q(n,d)$, $n \ge 2d$ \\
$d$ & $q$ & 0 & $q^e$ & dual polar graphs $C_q(d,e)$,
  $e=0,\frac12,1,\frac32,2$ \\
$d$ & $q$ & $q-1$ & $q^e-1$ & bilinear forms graph $H_q(d,e)$ \\
$\lfloor n/2 \rfloor$ & $q^2$ & $q^2-1$ & $q^{2n-2d-1}{-}1$ &
  alternating forms graphs $A_q(n)$ \\
$d$ & $-q$ & $-q-1$ & $-(-q)^d-1$ &
  Hermitian forms graphs $Q_q(d)$
\end{tabular}

\medskip
Below we give the asymptotic behavior of the eigenmatrix $P$
of these graphs when $d,b,\alpha$ are fixed and $\beta$ tends to infinity
(Theorem \ref{largebeta}).
We also give a simple explicit expression for the eigenvalues $P_{dj}$,
that perhaps has not been noticed before (Proposition \ref{lastrow}).

Subsequently, we investigate each of the individual families,
and determine smallest and second largest eigenvalues
and/or other properties of the eigenvalues.
Main results are Theorem \ref{Gr-smallest} for the Grassmann graphs,
Corollary \ref{dualpolar-smallest} for the dual polar graphs,
Theorem \ref{bilinear-2ndlargest} for the bilinear forms graphs,
Theorem \ref{alt-2ndlargest} for the alternating forms graphs,
and Theorem \ref{herm-smallest-2ndlargest} for the Hermitian forms graphs.

\section{The Hamming case}
We prove the stated results for the Hamming graphs.

\subsection{Identities}
We collect some (well-known) identities used in the sequel.

The defining equation gives $K_j(i)$ as a polynomial in $i$ of degree $j$
with leading coefficient $(-q)^j / j!$. We give three expressions.
\begin{align*}
K_j(i) &= \sum_{h=0}^j (-1)^h (q-1)^{j-h} \binom{i}{h} \binom{d-i}{j-h}\\
&= \sum_{h=0}^j (-q)^h (q-1)^{j-h} \binom{i}{h} \binom{d-h}{j-h}\\
&= \sum_{h=0}^j (-1)^h q^{j-h} \binom{d-i}{j-h} \binom{d-j+h}{h} 
\end{align*}
(see Delsarte \cite[p.~39]{Delsarte73}, and \cite[(15)]{Delsarte76}).

%
%
%
One has the symmetry
$$K_j(i) / \tbinom{d}{j} (q-1)^j = K_i(j) / \tbinom{d}{i} (q-1)^i.$$
In particular, $K_j(i)$ and $K_i(j)$ have the same sign.

There is also the symmetry
$$K_{d-j}(i) = (-1)^{i-j} (q-1)^{d-i-j} K_j(d-i).$$

\begin{Proposition} \label{prop:3term_hamming}
 Let $i,j \ge 1$. Then

$(q-1)(d-i)K_j(i+1) - (i+(q-1)(d-i)-qj)K_j(i) + iK_j(i-1) = 0$.
\end{Proposition}

%


%
%
%
%

\subsection{Proofs}
The occurrence of $d - \frac{d-1}{q}$ in Conjecture \ref{conj:vandam_sotirov} is explained by the following proposition.
Where it refers to $K_j(1)$ or $K_j(2)$, it is assumed that $d \ge 1$ or $d \ge 2$.

\begin{Proposition}\label{prop12}
Let $q \ge 2$ and $0 \le j \le d$.

(i) $K_j(1) < 0$ if and only if $j \ge d - \frac{d-1}{q}$.

(ii) $K_j(2) = K_j(1)$ if and only if $j = 0$ or $j = d - \frac{d-1}{q}$.

(ii)$'$ $K_j(2) > K_j(1)$ if and only if $j > d - \frac{d-1}{q}$.

(iii) $K_j(2) = \frac{-1}{q-1}K_j(1)$ if and only if $j = (d-1)(1 - \frac1q)$ or $j = d$.

(iv) Let $d - \frac{d-1}{q} \le j \le d$. Then $|K_j(2)| \le |K_j(1)|$.
\end{Proposition}

\Proof
(i)
Since $K_j(i)$ has the same sign as $K_i(j)$, this follows from
$K_1(j) = (q-1)d-qj$.

(ii)
Since $K_j(i) = \binom{d}{j} (q-1)^{j-i} K_i(j) / \binom{d}{i}$,
the claim says that
$K_2(j) = \frac12 (q-1)(d-1) K_1(j)$
precisely for the two specified values of $j$. But this condition
is quadratic in $j$, and is up to a constant factor $j(j-d+\frac{d-1}{q}) = 0$.

(ii)$'$
Clear from (ii), since $K_2(j)$ has positive leading coefficient.

(iii)
The condition is equivalent to $K_2(j) = -\frac12 (d-1) K_1(j)$.
Again it is quadratic in $j$. Up to a constant factor it is
$(j-d)(j-(d-1)(1 - \frac1q)) = 0$.

(iv)
We want to show that
$|K_2(j)| \le \frac{1}{2}(q-1)(d-1) |K_1(j)|$.
Since $K_1(j) < 0$ this is the pair of conditions
$K_2(j) - \frac{1}{2}(q-1)(d-1)K_1(j) \ge 0$ and
$-K_2(j) - \frac{1}{2}(q-1)(d-1)K_1(j) \ge 0$.

The former is up to a positive constant factor
equivalent to $j(j-d+\frac{d-1}{q}) \ge 0$.

For the latter it suffices to see that
$-K_2(j) - \frac{1}{2}(d-1)K_1(j) \ge 0$.
Up to a positive constant factor this is equivalent to
$(j-d)(j-(d-1)(1 - \frac1q)) \le 0$.
\qed

\medskip\noindent
If $j = d - \frac{d-1}{q}$, then $K_1(j) = -1$, and $K_j(1) = -\frac{1}{d} \binom{d}{j} (q-1)^{j-1}$.

\medskip
In order to prove Theorems \ref{binary} and \ref{nonbinary},
we need three lemmas.

\begin{Lemma}\label{qpow}
$|K_j(i)| \le (q-1)^{d-i} \binom{d}{j}$.
\end{Lemma}
\Proof
Since $\binom{d-i}{j-h} = 0$ unless $j-h \le d-i$, we have\\
$|K_j(i)| = |\sum_h (-1)^h (q-1)^{j-h} \binom{i}{h} \binom{d-i}{j-h}|
\le \sum_{h \ge i+j-d} (q-1)^{j-h} \binom{i}{h} \binom{d-i}{j-h}$\\
\phantom{$|K_j(i)|$} $\le (q-1)^{d-i} \sum_h \binom{i}{h} \binom{d-i}{j-h} =
(q-1)^{d-i} \binom{d}{j}$.
\qed

\begin{Lemma}\label{3term}
Let $1 < i < d$ and $d - \frac{d-1}{q} \le j \le d$.
If $qj \le 2(q-1)(d-i)$, then
$|K_j(i+1)| \le \max (|K_j(i-1)|,~|K_j(i)|)$.
\end{Lemma}
\Proof
Apply Proposition \ref{prop:3term_hamming}. Put $a = (q-1)(d-i)$.
One has $a K_j(i+1) - (i-qj+a)K_j(i) + iK_j(i-1) = 0$.
If $|K_j(i-1)| \le M$ and $|K_j(i)| \le M$, then
$a |K_j(i+1)| \le |i-qj+a| M + iM$, and the conclusion follows
if $i + |i-qj+a| \le a$.
Now $qj-i-a > (q-2)i \ge 0$, so we need $qj \le 2a$,
and that was one of the hypotheses.
\qed

\medskip
For $q=2$ the scheme is imprimitive, and the graphs $H(d,q,j)$
are bipartite for odd $j$, and disconnected for even $j$.
One has the additional symmetry $K_j(d-i) = (-i)^j K_j(i)$.

\begin{Lemma}\label{bd} Let $j < d/2$ and $0 < i < d$. Then
$\binom{d-1}{j-1} \le \sum_g \binom{i}{2g}\binom{d-i}{j-2g} \le \binom{d-1}{j}$.
\end{Lemma}

We prove Lemma \ref{bd} in the proof of Theorem \ref{binary}.

\begin{nTheorem}{\bf\ref{binary}} {\rm (\cite[Cor.~10]{DK})} Let $q = 2$.

(i) If $j \ne d/2$, then $|K_j(i)| \le |K_j(1)|$ for all $i$, $1 \le i \le d-1$.

(ii) If $j = d/2$, then $K_j(1) = 0$ and $|K_j(i)| \le |K_j(2)|$
for all $i$, $1 \le i \le d-1$.
\end{nTheorem}
\Proof
(i)
By the symmetry $K_{d-j}(i) = (-1)^i K_j(i)$ we may suppose $j < d/2$.

We prove Lemma \ref{bd} and part (i) of the theorem simultaneously. Since
$K_j(i) = \sum_h (-1)^h \binom{i}{h} \binom{d-i}{j-h}
= 2 \sum_g \binom{i}{2g} \binom{d-i}{j-2g} - \binom{d}{j}$,
and $K_j(1) = \binom{d-1}{j} - \binom{d-1}{j-1}$,
both statements are equivalent for all $i$.

Prove the statement of the lemma by induction of $d$.
The conclusion follows by adding the inequalities
for $(d-1,j-1)$ and $(d-1,j)$, using that $\binom{n}{m} = \binom{n-1}{m-1} + \binom{n-1}{m}$,
except possibly when $i = d-1$ or $j = (d-1)/2$.
If $i = d-1$, the claim is that
$\binom{d-1}{j-1} \le \binom{d-1}{2[j/2]} \le \binom{d-1}{j}$,
which is true.
Instead of treating $j = (d-1)/2$ we use symmetry and take $j = (d+1)/2$
and prove the statement in (i) by induction on $i$,
using Proposition \ref{prop12} (iv) and Lemma \ref{3term}.
Here we may suppose $2 \le i \le d/2$ by the symmetry $K_j(d-i) = (-1)^j K_j(i)$.

(ii)
By symmetry, $K_j(i) = 0$ when $j = d/2$ and $i$ is odd.
The 3-term recurrence reduces to $(d-i)K_j(i+1) + iK_j(i-1) = 0$
for odd $i$, so that $K_j(2h) = (-1)^h \binom{d}{d/2} \binom{d/2}{h} / \binom{d}{2h}$
and $|K_j(2h)|$ decreases with increasing \mbox{$2h \le d/2$}.
\qed

\begin{nCorollary}{\bf\ref{binarycor}} Let $q=2$ and $j \ge (d+1)/2$.

(i) One has $K_j(1) \le K_j(i)$ for all $i$, $0 \le i \le d-1$.

(ii) One has $K_j(1) \le K_j(d)$ if and only if $j$ is even or $j = d$.
\end{nCorollary}
\Proof
Since $K_j(1) < 0$, part (i) follows from part (i) of the theorem,
and part (ii) from $K_j(d) = (-1)^j K_j(0)$.
\qed

\medskip
Next, consider the nonbinary case.

\begin{nTheorem}{\bf\ref{nonbinary}}
Let $q \ge 3$ and $d - \frac{d-1}{q} \le j \le d$.

(i) One has $K_j(1) \le K_j(i)$ for all $i$, $0 \le i \le d$.

(ii) One has $|K_j(i)| \le |K_j(1)|$ for all $i \ge 1$,
unless $(q,d,i,j) = (3,4,3,3)$.
\end{nTheorem}

\noindent
If $(q,d,j) = (3,4,3)$ then $K_j(0) = 32$, $K_j(i) = -4$ for $i=1,2,4$,
and $K_j(3) = 5$.

\medskip\Proof
Since $K_j(1) < 0$ (and $K_j(0)$ is the largest eigenvalue), part (i)
follows from part (ii).
The case $i = 2$ was handled in Proposition \ref{prop12}, so we may assume $i \ge 3$.

For $j = d$ one has $K_j(i) = (-1)^i (q-1)^{d-i}$, and the statement is true.

For $j = d-1$ one has $K_j(i) = (-1)^{i-1} (q-1)^{d-i-1}(qi-d)$
and $d \ge q+1$.
To show the claim it suffices to show that $qi-d \le (q-1)^{i-1}(d-q)$ $(*)$,
and this follows from $q(i-1)-1 \le (q-1)^{i-1}$, unless $(q,i) = (3,3)$,
in which case $(*)$ still holds, unless $d = 4$.

So, we may assume $d-\frac{d-1}{q} \le j \le d-2$.
This implies that $3 \le q \le (d-1)/2$.

If $qj \le 2(q-1)(d-i+1)$ then we can apply Lemma \ref{3term}
(and induction on $i$) to conclude that $|K_j(i)| \le \max(|K_j(1)|,|K_j(2)|)$,
and we are done. So, assume $qj > 2(q-1)(d-i+1)$.

One has $K_j(1) = (q-1)^{j-1} \binom{d}{j} (q-1-\frac{qj}{d})$,
where the last factor is negative. From Lemma \ref{qpow} we see that
$|K_j(i)| \le |K_j(1)|$ when 
$d \le (q-1)^{i+j-d-1} (qj-(q-1)d)$.

Using $qj-(q-1)d \ge 1$
and $d-i+1 < \frac{qj}{2(q-1)} \le \frac34 j$
and $j \ge d - \frac{d-1}{q} \ge \frac23 d$ and $q \ge 3$
we see that it suffices to have $d^6 \le 2^d$, so $d \ge 30$ suffices.
The finitely many $d$ with $d < 30$ can be checked separately.
\qed

\subsection{Large \texorpdfstring{$q$}{q}}
\begin{Proposition} \label{hamming-large}
For fixed $d$, let $q$ be sufficiently large.
Then $K_j(i)$ is positive for $i+j \le d$, and has sign $(-1)^{i+j-d}$ for $i+j \ge d$.
For each $j > 0$, the smallest eigenvalue of $H(d,q,j)$ is $K_j(d-j+1)$.
\end{Proposition}
\Proof
We have $K_j(i) = \sum_{h=0}^j (-1)^h (q-1)^{j-h} \binom{i}{h} \binom{d-i}{j-h}$.
When $q$ tends to infinity, and $d,j$ are fixed, this sum is dominated by its
first nonzero term. So $K_j(i) \approx (q-1)^j \binom{d-i}{j}$ if $i+j \le d$,
and $K_j(i) \approx (-1)^{j+i-d} (q-1)^{d-i} \binom{i}{j+i-d}$ if $i+j \ge d$.
\qed

\medskip
How large is `sufficiently large'?
The value $K_j(d-j+1)$ is the unique smallest
eigenvalue of $H(d,q,j)$ for all $j$ when $q \ge q_0(d)$.

\medskip
\begin{tabular}{@{}c|c@{~~}c@{~~}c@{~~}c@{~~}c@{~~}c@{~~}c@{~~}c@{~~}c@{~~}c@{~~}c@{~~}c@{~~}c@{~~}c@{~~}c@{~~}c@{~~}c@{~~}c@{~}c@{}}
$d$ &   2 & 3 & 4 & 5 & 6 & 7 &  8 &  9 & 10 & 12 & 14 & 16 & 18 & 20 & 30 & 40 & 50 & 60 & 100 \\
\hline
$q_0$ & 2 & 3 & 4 & 5 & 7 & 9 & 12 & 15 & 18 & 26 & 35 & 45 & 57 & 70 & 156 & 277 & 433 & 623 & 1730
\end{tabular}

\begin{Lemma} Suppose $q > \frac14 d^2 + 1$. Then

(i) $K_j(i) > 0$ for $i \le d-j$,

(ii) $K_j(d-j+1) < 0$,

(iii) $|K_j(i)| < |K_j(d-j+1)|$ for $i > d-j+1$.
\end{Lemma}
\Proof
If $q > \frac14 d^2 + 1$, then the terms
$(q-1)^{j-h} \binom{i}{h} \binom{d-i}{j-h}$
decrease monotonically when $h$ increases, so that the sign of $K_j(i)$
is that of the first nonzero term and the difference between $K_j(i)$
and the first nonzero term is smaller than the next term.

For $2 \le e \le j$ we have\\
$|K_j(d-j+e)| \le (q-1)^{j-e} \binom{d-j+e}{e} + (q-1)^{j-e-1} \binom{d+j+e}{e+1} (j-e)$
and\\
$|K_j(d-j+1)| \ge (q-1)^{j-1} (d-j+1) - (q-1)^{j-2}\binom{d-j+1}{2} (j-1)$,
so that\\
$|K_j(d-j+e)|/(q-1)^{j-e-1} \le \binom{d-j+e}{e} (q-1+\frac{d-j}{e+1}(j-e))
\le \frac43 q \binom{d-j+e}{e}$ \\
and $|K_j(d-j+1)| \ge \frac12 q (q-1)^{j-2} (d-j+1)$.
So, it suffices to see\\
$\binom{d-j+e}{e} \le \frac38 (q-1)^{e-1} (d-j+1)$.
This holds for $e \ge 3$, and for $e=2$, $j \ge 3$, and for $j=e=2$ we can
drop the factor $\frac43$, and the conclusion holds.
\qed

\subsection{Coincidences}
A general matrix $A$ in the Bose-Mesner algebra ${\cal A}$ of a $d$-class
association scheme (see \cite[Chapter 2]{BCN} for a definition)
will have $d+1$ distinct eigenvalues, and generate ${\cal A}$,
in the sense that each element of ${\cal A}$ is a polynomial
of degree at most $d$ in $A$. Cases where some relation matrix $A_j$
has fewer eigenvalues (and hence generates a proper subalgebra)
are of interest.

Look at the Hamming scheme. For $q = 2$, the main expected coincidences
between the $P_{ij} = K_j(i)$ for fixed $d$ and $j$ are given
in the following lemma.

\begin{Lemma}\label{coin2} Let $q = 2$.

(i) If $j$ is even, then $P_{ij} = P_{d-i,j}$.

(ii) If $d = 2j$, then $P_{ij} = 0$ for all odd $i$.

(iii) If $d = 2j-1$, then $P_{2h-1,j} = P_{2h,j}$ for $1 \le h \le j-1$.

(iv) If $j = d$, then $P_{ij} = (-1)^i$ for all $i$.

\end{Lemma}

\Proof
We only have to show (iii), and this follows from
Proposition \ref{prop12} (ii), and the 3-term recurrence
given in Proposition \eqref{prop:3term_hamming}.
\qed

%
%
%

\medskip
If $K_j(i) = 0$, then also $K_j(d-i) = 0$ and we have a further
coincidence (when $j$ is odd and  $i \ne d/2$).
Integral zeros of Krawtchouk polynomials play a role e.g. in the study
of the existence of perfect codes or the invertibility of Radon transforms,
and have been studied by many authors,
cf.~\cite{CS,DG,HS,Hanrot,KL,SdW96,SdW99}.
For $j=1,2,3$ there are infinite families. For fixed $j \ge 4$ there are
zeros only for finitely many $d$. Recall that $K_j(i) = 0$ if and only if
$K_i(j) = 0$.

\begin{Lemma} {\rm (\cite[Th. 4.6]{CS} and \cite[Ex. 10]{DG})} Let $q=2$, $i \le d/2$, $j \le d/2$.

(i) $K_1(i) = 0$ if and only if $d = 2i$.

(ii) $K_2(i) = 0$ if and only if $i = \binom{h}{2}$, $d = h^2$ for some
integral $h \ge 3$.

(iii) $K_3(i) = 0$ if and only if $i = h(3h \pm 1)/2$,
$d = 3h^2+3h+\frac32 \pm (h+\frac12)$ for some integral $h \ge 2$.

(iv) $K_{2h}(4h-1) = 0$ if $d = 8h+1$.
\end{Lemma}

The family given last has $j = (d-3)/2$. There are also infinite families
with $j = (d-t)/2$ for $t=4,5,6,8$ (\cite{HS}).

\medskip
For arbitrary $q$ there are fewer obvious coincidences.

\begin{Lemma}\label{coinq} Let $q \ge 2$.

(i) If $j = 0$, then $P_{ij} = 1$ for all $i$.

(ii) If $j = 2$, then $P_{hj} = P_{ij}$ if and only if
$h+i = 2(d-1)(1-\frac{1}{q})+1$.

(iii) If $qj = (q-1)d+1$, then $P_{1j} = P_{2j}$.
\end{Lemma}

\Proof
(i) The matrix $A_0 = I$ only has the single eigenvalue 1.

(ii) Note that $K_2(i)$ is quadratic in $i$.

(iii) This is what Proposition \ref{prop12} (ii) says.
\qed

\medskip
We look for cases where some $A_j$ has fewer distinct eigenvalues
than expected (given the above lemmas), or just has few distinct eigenvalues.
Below we list cases where $H(d,q,j)$ has precisely $n$
distinct eigenvalues, while $d+1 > n$, for $n=3,4,5,6$.

\begin{Conjecture}
If $H(d,q,j)$ is connected, it has more than $d/2$ distinct eigenvalues.
\end{Conjecture}

\subsubsection{Three distinct eigenvalues}
If $H(d,q,j)$ has three distinct eigenvalues, it is strongly regular,
or (in case $q=2$ and $j$ even) it is the disjoint union of two
isomorphic connected components, both strongly regular.

For example, the $P$-matrix of $H(4,3)$ was given above,
$$
P = \left(\begin{array}{ccccc}
 1 &   8 &  24 & 32 & 16 \\
 1 &   5 &   6 & -4 & -8 \\
 1 &   2 &  -3 & -4 &  4 \\
 1 &  -1 &  -3 &  5 & -2 \\
 1 &  -4 &   6 & -4 &  1
\end{array}\right)
$$
and $H(4,3,3)$ is strongly regular with parameters
$(v,k,\lambda,\mu) = (81,32,13,12)$ and spectrum $32^1~5^{32}~(-4)^{48}$.

For $H(7,2)$ one gets

{\small 
$$
P = \left(\begin{array}{cccccccc}
1 &  7 & 21 &  35 & 35 &  21 &  7 &  1 \\
1 &  5 &  9 &   5 & -5 &  -9 & -5 & -1 \\
1 &  3 &  1 &  -5 & -5 &   1 &  3 &  1 \\
1 &  1 & -3 &  -3 &  3 &   3 & -1 & -1 \\
1 & -1 & -3 &   3 &  3 &  -3 & -1 &  1 \\
1 & -3 &  1 &   5 & -5 &  -1 &  3 & -1 \\
1 & -5 &  9 &  -5 & -5 &   9 & -5 &  1 \\
1 & -7 & 21 & -35 & 35 & -21 &  7 & -1
\end{array}\right)
$$

}

\noindent
and the graph $H(7,2,4)$ has two connected components, both isomorphic
to the graph $\Delta$ on the 64 binary vectors of length 7 and even weight,
adjacent when they differ in 4 places. The graph $\Delta$ is strongly regular
with parameters $(v,k,\lambda,\mu) = (64,35,18,20)$
and spectrum $35^1~3^{35}~(-5)^{28}$.

\medskip
Cases with three eigenvalues (the connected graphs among these are
strongly regular---we give the standard parameters $(v,k,\lambda,\mu)$):

\medskip
\begin{tabular}{cccl}
$d$ & $q$ & $j$ & comment \\
\hline
4 & 2 & 2 & 2 copies of $\overline{4K_2}$\rule{0pt}{2.4ex} \\
5 & 2 & 2 & 2 copies of the Clebsch graph \\
5 & 2 & 4 & 2 copies of the complement of the Clebsch graph \\
7 & 2 & 4 & 2 copies of $VO^+(6,2)$ \\
\hline
4 & 3 & 2 & $(81,24,9,6)$ \\
4 & 3 & 3 & $(81,32,13,12)$: $VO^+(4,3)$ \\
\hline
3 & 4 & 2 & $(64,27,10,12)$: $VO^-(6,2)$ \\
\hline
\end{tabular}

\medskip
More generally, if we take the Hamming scheme $H(d,q)$ with $q=4$,
and call two distinct vertices adjacent if their distance is even,
we obtain a strongly regular graph (as was observed in \cite[Case III]{KagSD}),
namely the graph $VO^\pm(2d,2)$, where the sign is $(-1)^d$.
Indeed, the weight of a quaternary digit is given by the (elliptic)
binary quadratic form $x_1^2 + x_1x_2 + x_2^2$.
For $d=3$ this graph is $H(3,4,2)$.

\subsubsection{Four/five/six distinct eigenvalues}
In Table \ref{table} below we list further cases in which $H(d,q,j)$
has fewer than $d+1$ distinct eigenvalues.

\begin{table}[ht]
\centering
\begin{minipage}{.30\textwidth}\footnotesize
\begin{tabular}{@{}ccc@{~~}l@{}}
$d$ & $q$ & $j$ & comment \\
\hline
5 & 2 & 3 & L\ref{coin2} (iii) \\
6 & 2 & 2,4 & L\ref{coin2} (i) \\
7 & 2 & 2,6 & L\ref{coin2} (i) \\
8 & 2 & 4 & L\ref{coin2} (i),(ii) \\
11 & 2 & 6 & L\ref{coin2} (i),(iii) \\
\hline
5 & 3 & 3 & $P_{13} = P_{43}$ \\
  &   &   & $P_{23} = P_{53}$ \\
\hline
5 & 4 & 2 & L\ref{coinq} (ii) \\
\hline
4 & 6 & 2 & L\ref{coinq} (ii) \\
\hline
\end{tabular}
\caption*{Four eigenvalues}
\end{minipage}\hfill
\begin{minipage}{.33\textwidth}\footnotesize
\begin{tabular}{@{}ccc@{~~}l@{}}
$d$ & $q$ & $j$ & comment \\
\hline
6 & 2 & 3 & L\ref{coin2} (ii) \\
8 & 2 & 2,6 & L\ref{coin2} (i) \\
9 & 2 & 2,4,6,8 & L\ref{coin2} (i) \\
10 & 2 & 4 & L\ref{coin2} (i) \\
   &   &   & $P_{24} = P_{34}$ \\
10 & 2 & 8 & L\ref{coin2} (i) \\
   &   &   & $P_{38} = P_{48}$ \\
11 & 2 & 4 & L\ref{coin2} (i) \\
   &   &   & $P_{24} = P_{44}$ \\
11 & 2 & 8 & L\ref{coin2} (i) \\
   &   &   & $P_{38} = P_{58}$ \\
12 & 2 & 6 & L\ref{coin2} (i),(ii) \\
15 & 2 & 8 & L\ref{coin2} (i),(iii) \\
\hline
7 & 3 & 2 & L\ref{coinq} (ii) \\
7 & 3 & 5 & L\ref{coinq} (iii) \\
   &   &   & $P_{35} = P_{65}$ \\
   &   &   & $P_{55} = P_{75}$ \\
\hline
5 & 4 & 3 & $P_{33} = P_{53}$ \\
5 & 4 & 4 & L\ref{coinq} (iii) \\
\hline
6 & 5 & 2 & L\ref{coinq} (ii) \\
\hline
5 & 6 & 3 & $P_{23} = P_{53}$ \\
\hline
5 & 8 & 2 & L\ref{coinq} (ii) \\
\hline
\end{tabular}
\caption*{Five eigenvalues}
\end{minipage}\hfill
\begin{minipage}{.36\textwidth}\footnotesize
\begin{tabular}{@{}ccc@{~~}l@{}}
$d$ & $q$ & $j$ & comment \\
\hline
7 & 2 & 3 & $P_{13} = P_{53}$ \\
9 & 2 & 5 & L\ref{coin2} (iii) \\
10 & 2 & 2,6 & L\ref{coin2} (i) \\
11 & 2 & 2,10 & L\ref{coin2} (i) \\
12 & 2 & 4 & L\ref{coin2} (i) \\
   &   &   & $P_{24} = P_{64}$ \\
12 & 2 & 8 & L\ref{coin2} (i) \\
   &   &   & $P_{28} = P_{68}$ \\
16 & 2 & 8 & L\ref{coin2} (i),(ii) \\
19 & 2 & 10 & L\ref{coin2} (i),(iii) \\
\hline
7 & 3 & 3 & $P_{23} = P_{53} = P_{63}$ \\
7 & 3 & 6 & $P_{26} = P_{36} = P_{56}$ \\
\hline
7 & 4 & 2 & L\ref{coinq} (ii) \\
7 & 4 & 4 & $P_{24} = P_{64}$ \\
  &   &   & $P_{54} = P_{74}$ \\
\hline
6 & 5 & 3 & $P_{43} = P_{63}$ \\
6 & 5 & 5 & L\ref{coinq} (iii) \\
\hline
7 & 6 & 2 & L\ref{coinq} (ii) \\
\hline
6 & 10 & 2 & L\ref{coinq} (ii) \\
6 & 10 & 3 & $P_{43} = P_{63}$ \\
\hline
\end{tabular}
\caption*{Six eigenvalues}
\end{minipage}\hfill
\caption{Cases where $H(d,q,j)$ has fewer than $d+1$ distinct eigenvalues}\label{table}
\end{table}

For example, the eigenmatrix of $H(7,3)$ is

{\small

$$
P = \left(\begin{array}{cccccccc}
1 &  14 & 84 & 280 & 560 & 672 & 448 & 128 \\
1 &  11 & 48 & 100 & 80 & -48 & -128 & -64 \\
1 &  8 &  21 & 10 & -40 & -48 & 16 & 32 \\
1 &  5 &  3 &  -17 & -16 & 24 & 16 & -16 \\
1 &  2 &  -6 & -8 & 17 & 6 &  -20 & 8 \\
1 &  -1 & -6 & 10 & 5 &  -21 & 16 & -4 \\
1 &  -4 & 3 &  10 & -25 & 24 & -11 & 2 \\
1 &  -7 & 21 & -35 & 35 & -21 & 7 &  -1
\end{array}\right) .
$$

}

\noindent
and we see coincidences in columns 2, 3, 5, 6.

\section{The Johnson case}
The eigenvalues of $J(n,d,j)$ are $P_{ij} = E_j(i)$ $(0 \le i,j \le d)$.
We give three expressions for the $E_j(i)$:
\begin{align*}
E_j(i) &= \sum_{h=0}^j  (-1)^h \binom{i}{h} \binom{d-i}{j-h} \binom{n-d-i}{j-h}
\\ &= \sum_{h=0}^j (-1)^{j-h} \binom{d-i}{h} \binom{d-h}{j-h} \binom{n-d-i+h}{h}
\\ &= \sum_{h=0}^i (-1)^{i-h} \binom{i}{h} \binom{d-h}{j} \binom{n-d-i+h}{n-d-j}
\end{align*}
(see Delsarte \cite[p.~48]{Delsarte73},
and Karloff \cite[Theorem 2.1]{Karloff}).

\subsection{Identities}

Using the second of the expressions given above for $E_j(i)$
we find the eigenvalues of the Kneser graph.

\begin{Proposition} {\rm (Lov\'asz \cite{Lovasz79})}
The eigenvalues of the Kneser graph are
$P_{id} = (-1)^i \binom{n-d-i}{d-i} = (-1)^i \binom{n-d-i}{n-2d}$.
\end{Proposition}
\Proof 
We use that $\binom{n+h}{h} = (-1)^h \binom{-n-1}{h}$
and $\sum_h \binom{a}{c-h} \binom{b}{h} = \binom{a+b}{c}$ and find
$P_{id} = \sum_h (-1)^{d-h} \binom{d-i}{h} \binom{n-d+h-i}{h} =
(-1)^d \sum_h \binom{d-i}{d-i-h} \binom{-n+d+i-1}{h}$\\
${}= (-1)^d \binom{-n+2d-1}{d-i} = (-1)^i \binom{n-d-i}{d-i}$.
\qed


\medskip
Let us write $E_j^{n,d}(i)$ instead of $E_j(i)$ when it is necessary
to make the dependence on $n$ and $d$ explicit. Now we have the
following induction.

\begin{Proposition}\label{Eberind}
Let $i,j \ge 1$. Then
$E_j^{n+2,d+1}(i) = E_j^{n,d}(i-1) - E_{j-1}^{n,d}(i-1)$.
\end{Proposition}
\Proof
Using $E^{n,d}_j(i) =
\sum_h (-1)^h \binom{i}{h} \binom{d-i}{j-h} \binom{n-d-i}{j-h}$
one sees that the claim reduces to $\binom{i}{h} = \binom{i-1}{h}+\binom{i-1}{h-1}$.
\qed

\medskip
There is a symmetry if $n = 2d$.

\begin{Proposition} \label{symmetry}
If $n = 2d$, then $E_{d-j}(i) = (-1)^i E_j(i)$.
In particular, if moreover $j = d/2$, $i$ odd, then $E_j(i) = 0$. \qed
\end{Proposition}

\subsection{Coincidences}

The association scheme on the set $X$ of partitions of a $2k$-set
into two $k$-sets has $\lfloor \frac12 k + 1 \rfloor$ relations $R_j$
(mutual intersection sizes $0+k$, $1+(k-1)$, ...,
$\lfloor \frac12 k \rfloor + \lceil \frac12 k \rceil$).
If one picks a fixed element in the $2k$-set, one sees that $(X,R_j)$
is isomorphic to the graph on the $(k-1)$-subsets of a $(2k-1)$-set,
adjacent when they meet in either $j-1$ or $k-j-1$ points.
Thus, in the Johnson scheme with $n = 2d+1$, the matrices
$A_j + A_{d-j+1}$ have not more than $(d+3)/2$ distinct eigenvalues.

\begin{Proposition} \label{2d+1-coinc}
Let $n = 2d+1$ and $j = (d+1)/2$ and $0 < t < d/2$.
Then $E^{n,d}_j(2t-1) = E^{n,d}_j(2t) = E_j^{n-1,d}(2t-1)$. \qed
\end{Proposition}

%

\subsection{Negative \texorpdfstring{$E_j(1)$}{Ej(1)}}
Let us write $e := n-d$ to make our formulas shorter and nicer.

\begin{Proposition} \label{johnsonevneg} Let $j > 0$. Then

(i) $E_j(1) = 0$ if and only if $j = de/n$.

(ii) $E_j(1) < 0$ if and only if $j > de/n$,

(iii) $E_j(1) = E_j(2)$ if and only if $j(n-1) = de$.

(iv) $E_j(1) < E_j(2)$ if and only if $j(n-1) > de$.
\end{Proposition}
\Proof
(i)-(ii)
We have $E_j(1) = (1 - \frac{jn}{de})
\binom{d}{j} \binom{e}{j}$.

(iii)-(iv)
Let $j > 1$. Writing out the expressions for $E_j(1)$ and $E_j(2)$,
dividing by $\binom{d-2}{j-2}\binom{e-2}{j-2}$,
multiplying by $j^2(j-1)^2$, and simplifying, we see that
$E_j(1) \le E_j(2)$ is equivalent to $j(n-1) \ge de$.
(There is a factor $n-2$, but $i = 2$ occurs only for $d \ge 2$, $n \ge 4$.)
For $j=1$ we have $E_1(i) = (d-i)(e-i) - i = de -i(n-i+1)$,
and $E_j(1) \le E_j(2)$ is equivalent to $n \le 2$, which is false.
\qed

\medskip
For $J(8,3)$ we have

$$
P = \left(\begin{array}{cccc}
1 & 15 & 30 & 10\\
1 & 7 & -2 & -6\\
1 & 1 & -5 & 3\\
1 & -3 & 3 & -1
\end{array}\right) .
$$

\subsection{Auxiliary results}
For any regular graph $\Gamma$ with adjacency matrix $A$,
the sum of the squares of the eigenvalues of $\Gamma$ (i.e., of $A$)
is the trace of $A^2$, which is $vk$, if $\Gamma$ has $v$ vertices
and is regular of valency $k$.
We apply this to $J(n,d,j)$, and find
$v k_j = \sum_{i=0}^d m_i E_j(i)^2$, where
$v = \binom{n}{d}$ is the number of vertices of $J(n,d)$,
$k_j = \binom{d}{j}\binom{e}{j}$ is the valency of $J(n,d,j)$
(with $e := n-d$),
and $m_i = \binom{n}{i}-\binom{n}{i-1}$ is the multiplicity
of the $i$-th eigenvalue (cf.~\cite[9.1.2]{BCN}).
It follows that $E_j(i)^2 \le vk_j / m_i$.

\medskip
We need to estimate $k_j$ close to its maximum value,
and use Chv\'atal's tail inequality for the hypergeometric distribution.

\begin{Lemma}\label{lem:concentrate_john}
  Let $I = ( \frac{de}{n} - \sqrt{d}, \frac{de}{n} + \sqrt{d} )$.
  Then $\sum_{j \in I} k_j \geq \frac{8}{11} v$.
\end{Lemma}
\Proof
  Consider the random variable $X$ that is $j$ with probability $k_j/v$.
  It has expected value $E(X) = \frac{de}{n}$.
  According to Chv\'atal \cite{Chvatal} (cf.~\cite{Skala}),
  $$Pr(|X - E(X)| \geq td) \leq 2 \exp(-2t^2d).$$
  Choosing $t = d^{-1/2}$ yields the assertion,
  as $1 - 2 \exp(-2) > \frac{8}{11}$.
%
%
%
\qed

\begin{Lemma} \label{Eji-bound}
Let $j_0 = \frac{de}{n}$, and let $j_0 \le j < j_0 + \frac{3}{2}$.
If $\frac{de}{n-1} \le j < d$ and $i \ge 3$, then $|E_j(i)| \le |E_j(1)|$.
\end{Lemma}
\Proof
We start with some observations that hold when $d$ is not too small.

(1) Since $\frac{de}{n-1} \le j \le d-1$, 
we find $de \le (d-1)(d+e-1)$, that is, $e \le (d-1)^2$.

(2) Since $n^3/d^2e^2$ decreases with $e$ for $e \le 2d$, and
increases with $e$ for $e \ge 2d$, it is maximal
for $e = (d-1)^2$ (for $d \ge 7$), so that
$n/j_0^2 \le (d^2-d+1)^3/d^2(d-1)^4 < 1 + \frac{3}{2d}$
(for $d \ge 10$).
%
%

(3) We show that $k_{j-1}/k_j < 3$ if $d \ge 10$.
Indeed, $k_{j-1}/k_j = c_j/b_{j-1} = j^2/(d-j+1)(e-j+1)$
so that $k_{j-1}/k_j < 3$ is equivalent to
$de-n(j-1)+(j-1)^2 - \frac13 j^2 > 0$.
The LHS decreases with $j$, so it suffices to show this
for $j = j_0 + \frac{3}{2}$. Since $de = nj_0$ we have to show
$\frac43 j_0^2 \ge n+1$, and this follows~from~(2).

(4) We show that $v / k_j < \frac16 (n-5)$ for $n \ge 42$.
%
%
According to Lemma \ref{lem:concentrate_john},
$\sum_{|\ell-j_0|<\sqrt{d}} k_{\ell} > \frac{8}{11} v$.
Let $k_{j_1}$ be the largest among the $k_{\ell}$.
Then $\lfloor j_0 \rfloor \le j_1 \le \lceil j_0 \rceil$
and $2\sqrt{d} \,\, k_{j_1} > \frac{8}{11} v$, that is,
$v / k_{j_1} <  \frac{11}{4} \sqrt{d}$.
The index $j$ differs at most 2 from $j_1$, and $j \ge j_1$.
Since $k_{j-2}/k_{j-1} \le k_{j-1}/k_j < 3$ we have
$k_{j_1} / k_j < 9$ and hence
$v / k_j < \frac{99}{4} \sqrt{d}$.
Our aim was $v / k_j < \frac16 (n-5)$, and since $n \ge 2d$
so that $n \ge \sqrt{2n} \sqrt{d}$ this follows from $n > 11255$.
The finitely many cases with $42 \le n \le 11255$ were checked by computer.
%
%

(5) We show that $E_j(i)^2 \le E_j(1)^2$ if $i \ge 3$.
In the discussion above we found that $E_j(i)^2 \le vk_j/m_i$,
where $m_i \ge m_3 = \frac16 n(n-1)(n-5)$. On the other hand,
$E_j(1) = \binom{d}{j} \binom{e}{j} (1 - \frac{j}{j_0})$,
and $j-j_0 \ge \frac{de}{n-1} - \frac{de}{n} = \frac{j_0}{n-1}$,
so that $E_j(i)^2 \le E_j(1)^2$ will hold when
$v/k_j \le \frac16 \frac{n(n-5)}{n-1}$. That was shown in (4).
Earlier we needed $d \ge 10$ (or $n \ge 42$), but if $d \le 9$ then
$n \le 73$, and these cases were checked by computer.
%
%
\qed

\subsection{The smallest eigenvalue}
It looks like $|E_j(1)|$ is the largest among the $|E_j(i)|$ ($1 \le i \le d$)
when $j$ is not very close to the zero $\frac{de}{n}$ of $E_j(1)$
(viewed as polynomial in $j$), say at least when
$|j - \frac{de}{n}| \ge \frac{1}{4}$.
If $|E_j(1)|$ is largest, and moreover $E_j(1) < 0$,
then $E_j(1)$ is the smallest among the $E_j(i)$, $0 \le i \le d$.
We prove below that this is the case if and only if $j \geq \frac{de}{n-1}$.

\begin{Lemma} \label{johnsonineq}
Let $(j-1)(n+1) \ge de$.
Then $E_j(0) + |E_{j-1}(1)| + |E_j(1)| \le E_{j-1}(0)$.
\end{Lemma}
\Proof
Use $E_j(0) = \binom{d}{j} \binom{e}{j}$ and
$E_j(1) = \binom{d}{j} \binom{e}{j} (1 - \frac{jn}{de})$
and $jn > de$ to see that the desired inequality is equivalent to
$\frac{jn}{de} \frac{d-j+1}{j} \frac{e-j+1}{j} +
|1 - \frac{(j-1)n}{de}| \le 1$.
If $(j-1)n \le de$ we have to show that
$(d-j+1)(e-j+1) \le j(j-1)$, that is, $de \le (j-1)(n+1)$,
which is our hypothesis.
If  $(j-1)n \ge de$ we have to show that
$(d-j+1) (e-j+1) + j(j-1) \le \frac{2dej}{n}$,
that is, $de - (j-1)(n-2j+1) \le \frac{2dej}{n}$,
that is, $de(n-2j) \le (j-1)n(n-2j+1)$, which holds by hypothesis. 
\qed

\medskip
Since we know the eigenvalues of the Kneser graph, the case $j = d$
is immediate.

\begin{Proposition} \label{johnson-d}
Let $d \ge 1$.
The smallest eigenvalue of $K(n,d)$, and the second largest
in absolute value, is $E_d(1)$. \qed
\end{Proposition}

\begin{Theorem} \label{smallestj1}
Let $j > 0$. Then $E_j(1)$ is the smallest eigenvalue of $J(n,d,j)$
if and only if $j(n-1) \ge de$. In this case $E_j(1)$ is also
the second largest in absolute value among the eigenvalues of $J(n,d,j)$.
\end{Theorem}
\Proof
By Proposition \ref{johnsonevneg}, if $E^{n,d}_j(1)$ is the smallest
eigenvalue of $J(n,d,j)$, then $j(n-1) \ge de$, and $E^{n,d}_j(1) < 0$.
We now show by induction on $d$ that if $j(n-1) \ge de$,
then $|E^{n,d}_j(i)| \le |E^{n,d}_j(1)|$. 
If $j = d$ the statement follows from Proposition \ref{johnson-d}.
If $\frac{de}{n-1} \le j < \frac{de}{n-3}$, then (since $n \ge 2d$
and $d \ge 3$) $j_0 < j < j_0 + \frac32$, where $j_0 = \frac{de}{n}$,
and our claim holds by Lemma \ref{Eji-bound} if $i \ge 3$.
We wish to~show that if $j(n-1) \ge (d+1)(e+1)$, then
$|E^{n+2,d+1}_j(i)| \le |E^{n+2,d+1}_j(1)|$, that is,
by Proposition \ref{Eberind},
$|E^{n,d}_j(i-1) - E^{n,d}_{j-1}(i-1)| \le |E^{n,d}_j(0) - E^{n,d}_{j-1}(0)|$.
Now $j(n-1) \ge (d+1)(e+1)$ implies $(j-1)(n-1) \ge de$, and
by induction, or trivially if $i=2$, $|E^{n,d}_j(i-1)| \le |E^{n,d}_j(1)|$ and
$|E^{n,d}_{j-1}(i-1)| \le E^{n,d}_{j-1}(1)|$ and our claim follows by
Lemma \ref{johnsonineq}.
\qed

\medskip
Karloff \cite{Karloff} studied graphs $J(n,d,j)$
for the special case $n = 2d$.
(His notation is $J(n,d,d-j)$ instead of our $J(n,d,j)$.)
He proves (\cite{Karloff}, Theorem 2.3) that $E_j(1)$
is the smallest eigenvalue of $J(n,d,j)$ when $d = n/2$ and $j \ge 5d/6$.
He conjectures (\cite{Karloff}, Conjecture 2.12) that $E_j(1)$
is the smallest eigenvalue of $J(n,d,j)$ when $d = n/2$ and $j > d/2$.
This conjecture immediately follows from the above theorem.

\begin{Corollary}\label{Kconj}
If $j > d/2$, then the smallest eigenvalue of $J(2d,d,j)$,
and the second largest in absolute value, is $E_j(1)$. \qed
\end{Corollary}

\medskip
For $n = 2d+1$ and $j = \frac12 d$ we have $E_j(2) = -\frac{d}{d-1} E_j(1)$
so that $|E_j(2)| > |E_j(1)|$.

\subsection{Large \texorpdfstring{$n$}{n}}
\begin{Proposition} \label{johnson-large}
For fixed $d$, let $n$ be sufficiently large.
Then $E_j(i)$ is positive for $i+j \le d$, and has sign $(-1)^{i+j-d}$ for $i+j \ge d$.
For each $j > 0$, the smallest eigenvalue of $J(n,d,j)$ is $E_j(d-j+1)$.
\end{Proposition}
\Proof
We have $E_j(i) =
\sum_{h=0}^j (-1)^h \binom{i}{h} \binom{d-i}{j-h} \binom{n-d-i}{j-h}$.
When $n$ tends to infinity, and $d$ is fixed, this sum is dominated by its
first nonzero term. So $E_j(i) \approx \binom{d-i}{j} \binom{n-d-i}{j}$
if $i+j \le d$, and
$E_j(i) \approx (-1)^{i+j-d} \binom{i}{i+j-d} \binom{n-d-i}{d-i}$
if $i+j \ge d$.
Also, for $i+j < d$,
$E_j(i)/E_j(i+1) \approx \frac{(d-i)(n-d-i)}{(d-i-j)(n-d-i-j)} > 1$,
so the $E_j(i)$ decrease in absolute value with increasing $i$.
\qed

\medskip
For example, for $J(27,5)$:

{\small

$$
P = \left(\begin{array}{cccccc}
1 & 110 & 2310 & 15400 & 36575 & 26334\\
1 & 83 & 1176 & 4060 & 665 & -5985\\
1 & 58 & 451 & 60 & -1710 & 1140\\
1 & 35 & 60 & -400 & 475 & -171\\
1 & 14 & -66 & 104 & -71 & 18\\
1 & -5 & 10 & -10 & 5 & -1
\end{array}\right) .
$$

}
For $d=5$ this is the smallest $n$ with the described sign pattern.
We have to go to $n=34$ to get decreasing absolute values in the columns.


%

\section{Graphs with classical parameters}\label{classparams}
Given a constant $b$, define
$$
\bnom{n}{m} = \bnom{n}{m}_b = \left\{
\begin{array}{ll}
\displaystyle 0 & \mbox{if $m < 0$,} \\
\displaystyle \binom{n}{m} & \mbox{if $b = 1$,} \\
\displaystyle \prod_{h=0}^{m-1} \frac{b^{n-h}-1}{b^{m-h}-1} & \mbox{otherwise.}
\end{array}\right.
$$

Graphs {\em with classical parameters} are
distance-regular graphs with intersection numbers
$b_i = (\bnom{d}{1} - \bnom{i}{1})(\beta - \alpha \bnom{i}{1})$ and
$c_i = \bnom{i}{1}(1 + \alpha\bnom{i-1}{1})$ ($0 \le i \le d$)
(see \cite[\S6.1]{BCN}).
It follows that $k = \beta \bnom{d}{1}$ and
$a_i = \bnom{i}{1} ( \beta-1+\alpha ( \bnom{d}{1} - \bnom{i}{1} - \bnom{i-1}{1}))$.
In \cite{BCN}, Corollary 8.4.2, the eigenvalues of graphs with
classical parameters are found to be
$\theta_i = \bnom{d-i}{1}(\beta - \alpha\bnom{i}{1}) - \bnom{i}{1}$
($0 \le i \le d$).

The base $b$ is an integer different from 0, $-1$ (\cite[6.2.1]{BCN}).

\subsection{Identities}
The $P_{ij}$ follow from the recurrence
$P_{i,j+1} = ((\theta_i - a_j) P_{ij} - b_{j-1} P_{i,j-1})/c_{j+1}$ and
the starting values $P_{i0} = 1$, $P_{i1} = \theta_i$ (see \cite[Chapter 4.1 (11)]{BCN}).
There is a simple explicit expression for the last row of the $P$ matrix.
It is independent of $\alpha$ and $\beta$.

\begin{Proposition}\label{lastrow}
$P_{dj} = (-1)^j \bnom{d}{j} b^{\binom{j}{2}}$.
\end{Proposition}
\Proof
Induction on $j$, using the recurrence.
\qed

\medskip
Graphs with classical parameters are formally self-dual
when $\alpha = b-1$. If this is the case, then
$P_{ij}/P_{0j} = P_{ji}/P_{0i}$ for all $i,j$,
and the number of vertices is $v = (\beta+1)^d$.
In this case, the above proposition can be translated to give
the values of the last column of $P$.

\begin{Proposition}\label{lastcol}
$P_{id} / P_{i+1,d} = 1 - (\beta+1) b^{-i}$. \qed
\end{Proposition}

\subsection{Sign changes}
The columns of the matrix $P$ correspond to the graph distances
on the distance-regular graph under consideration, and hence
have a natural ordering. For general distance-regular graphs
one is free to choose the ordering of the rows, corresponding
to an ordering of the eigenspaces. According to
\cite{BH}, Proposition 11.6.2, the $i$-th row and the $i$-th column
of $P$ have exactly $i$ sign changes if we order the rows
according to descending real order on the $\theta_i$.

Graphs with classical parameters are Q-polynomial, and hence
have a natural ordering on the eigenspaces. Usually this is
the order with descending $\theta_i$, provided $b > 0$.

\begin{Proposition} \label{sign-changes}
Suppose $b > 0$.
Then $\theta_0 > \theta_1 > \ldots > \theta_d$ if and only if
$\alpha \le b-1$ or $\beta > \alpha \bnom{d-1}{1} - b^{d-1}$.
If this is the case, then the $i$-th row and the $i$-th column of $P$
have exactly $i$ sign changes ($0 \le i \le d$).
\end{Proposition}
\Proof
We have to check that $\theta_i > \theta_{i+1}$, i.e., that
$\beta > \alpha \bnom{2i+1-d}{1} - b^{2i+1-d}$ for $0 \le i \le d-1$.
If $\alpha \le b-1$ then the strongest of these is the inequality
for $i=0$, but it is automatically satisfied since $\theta_0$
is the graph valency. If $\alpha > b-1$ the strongest is the
inequality for $i = d-1$, and we find the stated bound on $\beta$.
\qed

\medskip
The hypothesis of this proposition is satisfied for all families
of graphs with classical parameters considered in this note,
except for that of the Hermitian forms graphs, which have $b < 0$.

\medskip
In many cases the sign pattern is forced.

\begin{Proposition} \label{sign-pattern}
If the $i$-th row and the $i$-th column of $P$ have exactly $i$ sign changes,
and $P_{ij} > 0$ if $i+j \le d$, then $P_{ij}$ has sign
$(-1)^{i+j-d}$ if $i+j \ge d$.
\end{Proposition}
\Proof The only way to have $i$ sign changes in $P_{ij}$,
$d-i \le j \le d$ is to have $P_{ij}$ and $P_{i,j+1}$ of
opposite sign for all $j$, $d-i \le j \le d-1$. \qed

\subsection{Large \texorpdfstring{$\beta$}{beta}}
In the theorem below we show for graphs with classical parameters
$(d,b,\alpha,\beta)$ that if $(d,b,\alpha)$ is fixed and $\beta$
is large, then $P_{d-j+1,j}$ is the smallest eigenvalue of
the distance-$j$ graph, and $|P_{1j}|$ is its second largest
eigenvalue in absolute value.
We also determine the sign pattern of the matrix $P$.
This generalizes Propositions \ref{hamming-large} and
\ref{johnson-large} above.

There are families of graphs with classical parameters with $b < 1$,
such as the Hermitian forms graphs and the triality graphs.
However, Metsch \cite{Metsch99} showed that $\beta$ is bounded as a
function of $(d,b,\alpha)$ unless the graph is a Hamming, Johnson,
Grassmann, or bilinear forms graph. It follows that $b \ge 1$
when $\beta$ is unbounded.

\begin{Theorem} \label{largebeta}
For fixed $(d,b,\alpha)$, let $\beta$ be sufficiently large. Then

(i) $P_{ij} > 0$ for $i+j \le d$, and $P_{ij}$ has sign $(-1)^{i+j-d}$
for $i+j \ge d$.

(ii) $P_{d-j+1,j} = \min \{ P_{ij} \mid 0 \le i \le d \}$ for $j > 0$.

(iii) If $b \ge 1$, then $|P_{i+1,j}| < |P_{ij}|$ for $0 \le i \le d-1$.

\end{Theorem}
\Proof
For $|\beta| \to \infty$, we have $a_i \sim \bnom{i}{1}\beta$,
hence $\beta > 0$ and $b+1 \ge 0$ since $a_i \ge 0$ for $i=1,2$.
By \cite{BCN} (6.2.1), $b$ is an integer different from $0,-1$,
so $b \ge 1$.

(i)
In order to prove this, one only has to prove the first part,
then the second part follows by Propositions \ref{sign-changes}
and \ref{sign-pattern}.

From the recurrence $P_{i,j+1} = ((\theta_i - a_j) P_{ij} -
b_{j-1} P_{i,j-1})/c_{j+1}$ and
$b_i \sim (\bnom{d}{1} - \bnom{i}{1})\beta$, and
$c_i = O(1)$, and $a_i \sim \bnom{i}{1}\beta$, and
$\theta_i \sim \bnom{d-i}{1}\beta$,
it follows by induction that $P_{ij} \sim C_{ij} \beta^j$
for $i+j \le d$ and some positive constants $C_{ij}$.
%
%

(ii)
Now we know that $P_{d-j+1,j} < 0$ for large $\beta$.
By downward induction on $j$ one sees
that all $P_{ij}$ with $j \ge d-i$ have the same degree $m_i$ in $\beta$.
(Indeed, let $P_{id}$ have degree $m = m_i$ in $\beta$. Then
$c_{d+1-h}P_{i,d+1-h} = (\theta_i-a_{d-h})P_{i,d-h} - b_{d-h-1}P_{i,d-h-1}$
applied for $h=0,1,...,i-1$ shows that $P_{i,d-h-1}$ has degree $m$ in $\beta$
since the LHS has degree (at most) $m$, the middle term precisely $m+1$
and the final term must cancel that highest term.)
Since $P_{i,d-i}$ has degree $d-i$ this proves that $m_i = d-i$.
It follows that $P_{ij} \sim D_{ij} \beta^{d-i}$ for $i+j \ge d$
and some nonzero constants $D_{ij}$.
%
%
Thus, $P_{d-j+1,j}$ is the most negative in its column
when $\beta$ is large enough.

(iii)
In the interval $d-j \le i \le d$ the $P_{ij}$ have decreasing degrees
$d-i$ in $\beta$ and hence decrease in absolute value when $\beta$ is
sufficiently large. For the interval $0 \le i \le d-j$ the degree
is always $j$, and we have to work a bit more.

Put (just here) $c_{d+1} = 1$. Define polynomials
$F_j(x)$ for $-1 \le j \le d+1$ by $F_{-1}(x) = 0$, $F_0(x) = 1$,
$c_{j+1}F_{j+1}(x) = (x-a_j)F_j(x) - F_{j-1}(x) b_{j-1}$.
Then each $F_j$ has degree $j$ in $x$ (for $j \ge 0$),
and $P_{ij} = F_j(\theta_i)$ $(0 \le i,j \le d)$.
Finally, $F_{d+1}(\theta_i) = 0$ $(0 \le i \le d)$.
The $c_j$ are independent of $\beta$, but $a_j$ and $b_j$ and $\theta_i$
depend linearly on $\beta$.
Consider the coefficient of $\beta$ in
$\theta_i = \bnom{d-i}{1}(\beta - \alpha\bnom{i}{1}) - \bnom{i}{1}$
a linear expression in the variable $w = b^{-i}$ (if $b \ne 1$)
or $i$ (if $b = 1$).
Then the coefficient of $\beta^j$ in $F_j(\theta_i)$ is a degree $j$
polynomial $g_j(w) = \prod_{h=0}^{j-1} (\bnom{d-i}{1}-\bnom{h}{1})$
that vanishes for $d-j+1 \le i \le d$ and hence nowhere else.
That means that $P_{ij} = F_j(\theta_i) \sim g_j(b^{-i}) \beta^j$
(or $g_j(i)\beta^j$) is monotone in $i$, assuming $b \ge 1$.
Since $\bnom{d-i}{1}$ decreases with increasing $i$, also $P_{ij}$ does
(for $0 \le i \le d-j$).
\qed

\section{Grassmann graphs}
The Grassmann graphs $G_q(n,d)$ are the graphs with as vertices
the $d$-subspaces of an $n$-dimensional vector space over $\Ff_q$,
adjacent when they meet in codimension 1.
W.l.o.g. we assume $n \ge 2d$ (since $G_q(n,d)$ is isomorphic to $G_q(n,n-d)$),
and then these graphs are distance-regular of diameter $d$.
Let $G_q(n,d,j)$ be the distance-$j$ graph of $G_q(n,d)$,
where $0 \le j \le d$. The eigenvalues of $G_q(n,d,j)$ are
$P_{ij} = G_j(i)$ ($0 \le i \le d$), where
\enlargethispage*{0.3cm}
\begin{align*}
G_j(i) &= \sum_{h=0}^j (-1)^{j-h} \, q^{hi + \binom{j-h}{2}} \,
\qnom{d-i}{h} \qnom{d-h}{j-h} \qnom{n-d-i+h}{h} \\
&= \sum_{h=0}^i (-1)^{i-h} \, q^{j(j-i+h) + \binom{i-h}{2}} \,
\qnom{i}{h} \qnom{d-h}{j} \qnom{n-d-i+h}{n-d-j}
\end{align*}
(see Delsarte \cite{Delsarte76a}, Theorem 10,
and Eisfeld \cite{JE}, Theorem 2.7).

\subsection{Identities}

\begin{Proposition}
$G_d(i) = (-1)^i q^{d(d-i)+\binom{i}{2}} \qnom{n-d-i}{d-i}$. \qed
\end{Proposition}

Let us write $G_j^{n,d}(i)$ instead of $G_j(i)$ when it is necessary
to make the dependence on $n$ and $d$ explicit.
The analog of Proposition \ref{Eberind} is as follows.

\begin{Proposition}\label{Grassind}
Let $i,j \ge 1$. Then
$$G_j^{n+2,d+1}(i) = q^j G_j^{n,d}(i-1) - q^{j-1} G_{j-1}^{n,d}(i-1).$$
\end{Proposition}
\Proof
Use the first formula for $G^{n,d}_j(i)$, and
$\qnom{n+1}{m} = q^m \qnom{n}{m} + \qnom{n}{m-1}$.
\qed

\subsection{The smallest eigenvalue}
In Theorem \ref{Gr-smallest} we find the smallest among
the eigenvalues of $G_q(n,d,j)$ (for $(n,q) \ne (2d,2)$).
In Proposition \ref{prop:largest_absolute_grassmann} (ii)
we determine the second largest in absolute value (in all cases).

\medskip
The following lemma provides tools to estimate
Gaussian coefficients, and their quotients.

\begin{Lemma}\quad \label{gauss_bnds}

(i) If $n \le m$, $b > 1$, then $(b^n-1)/(b^m-1) \le b^{n-m}$.

(ii) If $m \ge 1$, $b > 1$, then
$(b^n-1)/(b^m-1) < b^{n-m+1}/(b-1)$.

(iii) If $b > 1$, then $\bnom{n}{k}_b \ge b^{k(n-k)}$.

(iv) {\rm (\cite[Lemma 37]{FIKM})}
If $0 < k < n$, $b > 1$, then $\bnom{n}{k}_b \ge (1+\frac{1}{b}) b^{k(n-k)}$.

(v) {\rm (\cite[Lemma 34]{FIKM})}
If $0 \le k \le n$, $b \ge 4$, then
$\bnom{n}{k}_b < (1+\frac{2}{b}) b^{k(n-k)}$. \qed
\end{Lemma}

\begin{Proposition}\quad \label{prop:largest_absolute_grassmann}

(i)
$G_j(1) < 0$ if and only if $j = d$.
$G_j(1)$ is never zero.

(ii) Let $i \ge 1$. Then $|G_j(i)| \le |G_j(1)|$.

(iii) Let $j \ge 1$, $i+j \le d$.
Then $0 < G_{j-1}(i) < G_j(i)$ if not $q=2$, $n=2d$, $i+j=d$.

(iv) Let $(n,q) \ne (2d,2)$.
Then $G_j(i)$ has sign $(-1)^{\max(0,i+j-d)}$.

(v) Among the $G_d(i)$ with $i \ge 0$, the smallest is $G_d(1)$.
\end{Proposition}
\Proof

(i) This is immediate from the second expression for $G_j(i)$.

(ii) Using $G_j(0) = q^{j^2} \qnom{d}{j} \qnom{e}{j}$ and
$G_j(1) = q^{j^2} \qnom{d-1}{j} \qnom{e}{j}
- q^{j(j-1)} \qnom{d}{j} \qnom{e-1}{j-1}$,
where $e = n-d$, we see that
$G_{j-1}(0) + |G_{j-1}(1)| + |qG_j(1)| \le qG_j(0)$.

Now apply induction on $d$ and $i$:
$|G^{n+2,d+1}(i)| \le |G^{n+2,d+1}_j(1)|$ follows from
$q|G^{n,d}_j(i-1)|+|G^{n,d}_{j-1}(i-1)| \le
q|G^{n,d}_j(1)|+|G^{n,d}_{j-1}(1)| \le qG^{n,d}_j(0)-G^{n,d}_{j-1}(0)$.

(iii) Induction on $d$. Positiveness follows from monotony since $G_0(i) = 1$.
For $i=0$ we have to show that $q^{j^2} \qnom{d}{j} \qnom{e}{j}$
increases with $j$, and it does, with the indicated exception.
Now for $i > 0$, using $j+1 \le d$ and $q \ge 2$:\\
$G_{j+1}^{n+2,d+1}(i)-G_j^{n+2,d+1}(i) =
q^{j+1}G^{n,d}_{j+1}(i-1)-2q^jG^{n,d}_j(i-1)+q^{j-1}G^{n,d}_{j-1}(i-1) > 0$.

(iv) This follows by part (iii) and Propositions \ref{sign-changes},
\ref{sign-pattern}.

(v) This follows by parts (i) and (ii).
\qed

\begin{Conjecture}\quad

(i) If $(n,q) \ne (2d,2)$, then $|G_j(i+1)| < |G_j(i)|$ when $0 \le i \le d-1$.

(ii) If $(n,q) = (2d,2)$, then $G_j(d-j)$ is negative
for $(d,j) = (5,3)$ and when $d \ge 6$, $2 \le j \le d-2$,
and $G_j(d-j)$ is the smallest among the $G_j(i)$ when
$d \ge 6$, $3 \le j \le d-2$.
\end{Conjecture}
%
%

We can prove part (i) for $q \ge 5$, but omit the details.

\medskip
We show that $G_j(i)$ is well-approximated by its main term $T$.

\begin{Lemma} \label{Gr-bound}
If $i+j \le d$, let $s := 1$ and
$T := q^{j^2} \qnom{d-i}{j} \qnom{n-d}{n-d-j}$.
If $i+j \ge d$, let $s := i+j-d$ and
$T := (-1)^s q^{j(d-i)+\binom{s}{2}} \qnom{i}{d-j} \qnom{n-i-j}{n-d-j}$.
If $q \ge 3$ or $q=2$, $n>2d$, then
$$
\left| \frac{G_j(i)}{T} - 1 \right| <
\frac{q^{2d+1-n}}{(q-1)^2} .
$$
\end{Lemma}
\Proof
Let $T_h$ be the term with index $h$ in the second expression for $G_j(i)$,
so that $T = T_m$ with $m = \min(i,d-j)$, and $0 \le h \le m$.
This expression is alternating, and
\begin{align*}
\left| \frac{T_{h-1}}{T_h} \right| &=
\left| -q^{-h+i-j} \frac{q^h-1}{q^{i-h+1}-1}
\frac{q^{d-h+1}-1}{q^{d-h-j+1}-1} \frac{q^{j-i+h}-1}{q^{n-d-i+h}-1}  \right| \\
  &< \frac{q^{d+h+j+1-n}}{(q-1)^2} \le \frac{q^{2d+1-n}}{(q-1)^2}
\end{align*}
if $h \ge 1$. (Here we used Lemma \ref{gauss_bnds} (ii) twice, and (i) once,
using that $h \le n-d-i+h$.)
If $q \ge 3$ or $q=2$, $n > 2d$, then the right-hand side is less than 1,
and the sum is alternating with decreasing terms, so that the difference
between the main term and the sum is not larger than the second term.
The main term is $T = T_m$, the maximal index that occurs.
\qed

\medskip
\Remark
For $q = 2$, $i \ge d-j+1$ we shall need a slightly
sharper bound. Now $i-h+1 \ge 2$ and in the inequalities
in the proof and conclusion of the lemma we can bound by
$q^{2d+2-n} / ((q-1)(q^2-1))$.

\medskip
Above the main term of the second expression for $G_j(i)$
was $T_i$ (if $i+j \le d$) or $T_{d-j}$ (if $i+j \ge d$).
If $q = 2$, $n = 2d$, $i+j \ge d \ge 6$, and $3 \le j \le d-2$,
the main term is $T_{d-j-1}$.


\begin{Lemma} \label{Gr-bound-sp-case}
Let $n=2d$, $q=2$, $d \geq 13$, $5 \leq j \leq d-5$ and $d-j \leq i < d$.
Set $s := i+j-d+1$.
Let $T := (-1)^s q^{j(d-i-1)+\binom{s}{2}} \qnom{i}{d-j-1} \qnom{j+1}{1} \qnom{2d-i-j-1}{d-j}$.
Then $|G_j(i)| \leq \frac{3}{2} |T|$.
For $i = d-j$, $G_j(i)$ is negative, and $|G_j(i)| \geq 5|T|/171$.
\end{Lemma}
\Proof
  Let $T_h$ be the term with index $h$ in the second expression for $G_j(i)$,
  so that $T = T_{d-j-1}$ and $0 \le h \le \min(i,d-j)$.
  As in the proof of Lemma \ref{Gr-bound}, we have
$$
\left| \frac{T_{h-1}}{T_h} \right| =
\left| -2^{-h+i-j} \frac{2^h-1}{2^{i-h+1}-1}
\frac{2^{d-h+1}-1}{2^{d-h-j+1}-1} \frac{2^{j-i+h}-1}{2^{d-i+h}-1}  \right|.
$$
  For $h \leq d-j-1$ ($\leq i-1$), we find using Lemma \ref{gauss_bnds} (i) with $h \le d-i+h$,
  \begin{align*}
    \left| \frac{T_{h-1}}{T_h} \right| &\leq 2^{2i-h-j-d} \cdot \frac{2^{j+2}}{3} \cdot \frac{2^{j-2i+2h+1}}{3}
    = \frac{2^{h +j+3-d}}{9} \leq \frac{4}{9}.
  \end{align*}
  For $h = d-j$ and $i+j > d$ we find, using $i \le d-1$ and $5 \le j \le d-5$,
  \begin{align*}
    \left| \frac{T}{T_{d-j}} \right| = 2^{i-d} \frac{(2^{d-j}-1)(2^{j+1}-1)(2^{d-i}-1)}{(2^{i+j-d+1}-1)(2^{2d-i-j}-1)} \ge \frac{31}{63}.
%
%
  \end{align*}
  For $h = d-j$ and $i+j= d$, we find, using $5 \leq j \leq d-5$
  and $d \geq 13$,
  \begin{align*}
    \left| \frac{T}{T_{d-j}} \right| &=
    2^{-j} \frac{(2^{d-j}-1)(2^{j+1}-1)(2^j-1)}{2^d-1}
    \ge \frac{31 \cdot 63 \cdot 255}{32 \cdot 8191} > \frac{19}{10}.
  \end{align*}
  If $i+j = d$, then $G_j(i) = \sum_{h=0}^{d-j} T_h$ is an alternating sum with terms
  increasing in absolute value up to $T = T_{d-j-1}$, and then decreasing again,
  hence $|\frac{G_j(i)}{T} - 1| \le \frac{4}{9} + \frac{10}{19} = 1 - \frac{5}{171} < 1$,
  so that $G_j(d-j)$ has the sign of $T$, i.e., is negative.
  For general $i$, if $G_j(i)$ has the same sign as $T$, then $|G_j(i)| \le |T|$.
  If $G_j(i)$ has the opposite sign, then
  $|G_j(i)| \le |T_{d-j}| - |T| + |T_{d-j-2}| \le
  (\frac{63}{31} - 1 + \frac{4}{9}) |T| < \frac{3}{2} |T|$.
\qed

\begin{Theorem} \label{Gr-smallest}
Let $1 \le j \le d$. 

(i) If $q \geq 3$ or $q=2$, $n \ge 2d+1$, then
the smallest eigenvalue of $G_q(n,d,j)$ is $G_j(d-j+1)$.

(ii) If $(n, q) = (2d, 2)$, and $7 \le j \le d-5$, then
the smallest eigenvalue of $G_q(n,d,j)$ is $G_j(d-j)$.
\end{Theorem}
\Proof
(i) The case $j = d$ is handled in
Proposition \ref{prop:largest_absolute_grassmann},
so we may assume $j < d$.
The smallest among the $G_j(i)$ is negative, and hence $i$ is one
of the values $d-j+1+2t$ where $t \ge 0$. First consider the case $q \ge 3$.
We compare $G_j(i)$ with $G_j(i+2)$.
By Lemma \ref{Gr-bound}, both are approximated by their main term $T$
with an error that is not larger than $\frac{3}{4}T$.
Let $T$, $T'$, $T''$ be the main terms for $G_j(i), G_j(i+1), G_j(i+2)$.
Then $|G_j(i+2)|/|G_j(i)| \le (\frac{7}{4}|T''|)/(\frac{1}{4}|T|)
= 7|T''|/|T|$. Now
$$
\frac{|T'|}{|T|} = q^{i-d} \frac{q^{i+1}-1}{q^{i-d+j+1}}
\frac{q^{d-i}-1}{q^{n-i-j}-1} <
\frac{q^{d+i-n+1}}{q-1}
$$
using Lemma \ref{gauss_bnds} (i), (ii), since $d-i \le n-i-j$.
It follows that $|T''|/|T| < (q^{2d+2i-2n+3})/(q-1)^2$.
Since $i+2 \le d$ and $n \ge 2d$ we have $2d+2i-2n+3 \le -1$
and $|G_j(i+2)|/|G_j(i)| \le 7|T''|/|T| < 7/12 < 1$, as desired.

For $q=2$, $n \ge 2d+1$ we use the remark following Lemma \ref{Gr-bound}
and find $|G_j(i+2)| \le \frac{5}{3}|T''|$ and
$|G_j(i)| \ge \frac{1}{3}|T|$, so that
$|G_j(i+2)|/|G_j(i)| \le 5|T''|/|T| < 5/8 < 1$, as desired.

(ii)
The cases with $d < 13$ can be checked by computer,
so we may assume $d \geq 13$.
The smallest among the $G_j(i)$ is negative, so has $i \ge d-j$
by Proposition \ref{prop:largest_absolute_grassmann} (iii).
The value $G_j(d-j)$ is negative. We show that it has maximal
absolute value among the $G_j(i)$ with $i \ge d-j$.

Let $T$ and $T'$ be the main terms of $G_j(d-j)$
and $G_j(i)$, where $i < d$. By Lemma \ref{Gr-bound-sp-case},
$|G_j(d-j)| \geq \frac{5}{171} |T|$
and $G_j(i) \leq \frac{3}{2} |T'|$.
Then $|G_j(i)/G_j(d-j)| \leq \frac{3 \cdot 171}{2 \cdot 5} |T'|/|T|$.
Now, as $d \geq 5$, for $i = d-j+1$ we have
\begin{align*}
  \frac{|T'|}{|T|} &= 2^{-j+1} \frac{\qnom{d-j+1}{2} \qnom{d-2}{j-2}}{\qnom{d-j}{1} \qnom{d-1}{j-1}}
  = 2^{-j+1} \frac{2^{d-j+1}-1}{3} \frac{2^{j-1}-1}{2^{d-1}-1} < \frac{2^{-j+2}}{3} \cdot \frac{16}{15}.
\end{align*}
As $\frac{2^{-j+2}}{3} \cdot \frac{16}{15} < \frac{2 \cdot 5}{3\cdot 171}$
for $j \ge 7$, $|G_j(d-j+1)| < |G_j(d-j)|$.
Now, let $d-j+1 \leq i \leq d-2$. Let $T'$ and $T''$ be the main terms
of $G_j(i)$ and $G_j(i+1)$. Then
\begin{align*}
  \frac{|T''|}{|T'|} &= 2^{i-d+1} \frac{2^{i+1}-1}{2^{i+j-d+2}-1}
                                  \frac{2^{d-i-1}-1}{2^{2d-i-j-1}-1}
  \leq \frac{4 \cdot 64}{3 \cdot 63} \, 2^{i-d} < 1.
\end{align*}
Hence, $|G_j(i)| < |G_j(d-j)|$ for $d-j+2 \leq i \leq d-1$.
Lemma \ref{Gr-bound-sp-case} excludes $i = d$, so we need to treat
that case separately. By Proposition \ref{lastrow},
$G_j(d) = (-1)^j \qnom{d}{j} q^{\binom{j}{2}}$, and hence
\begin{align*}
  \frac{|G_j(d)|}{|T|} &= \frac{2^{\binom{j}{2}} \, \qnom{d}{j}}
        {2^{j(j-1)} \, \qnom{d-j}{1}\qnom{j+1}{1}\qnom{d-1}{j-1}}
  \leq  \frac{2^{-\binom{j}{2}}}{(2^{j+1}-1)(2^j-1)} < 2^{-j}.
\end{align*}
Hence, $|G_j(d)| < |G_j(d-j)|$.
\qed

\section{Dual polar graphs}
Let $q$ be a prime power.
There are six types of finite classical polar spaces,
$C_d(q)$, $B_d(q)$, $D_d(q)$, $^2D_{d+1}(q)$, $^2A_{2d}(q)$,
and $^2A_{2d-1}(q)$ with associated parameter (in the same order)
$e = 1, 1, 0, 2, 1/2, 3/2$ (see \cite[\S9.4]{BCN}).
In the cases $^2A_{2d}(q)$ and $^2A_{2d-1}(q)$
the parameter $q$ is the square of a prime power.
The dual polar graphs $C_q(d, e)$ are the graphs with as vertices the
maximal subspaces of a polar space of rank $d$ with parameter $e$
over $\Ff_q$, adjacent when they meet in codimension 1.
These graphs are distance-regular of diameter $d$.
The eigenmatrix $P$ has entries $P_{ij} = C_j(i)$, where
\begin{align*}
  C_j(i) &= \sum_{h=\max(i-j, 0)}^{\min(d-j,i)} (-1)^{i-h} \,
         q^{\binom{i-h}{2} + \binom{j-i+h}{2} + (j-i+h) e} \,
         \qnom{d-i}{d-j-h} \qnom{i}{h} .
\end{align*}
This formula was taken from Vanhove \cite[Theorem 4.3.6]{FV}.
An expression in terms of $q$-Krawtchouk polynomials was
given in Stanton \cite[Thm. 5.4]{Stanton80}.

\subsection{Identities}

Let us write $C_j^{d}(i)$ instead of $C_j(i)$ when it is necessary
to make the dependence on $d$ explicit.

\begin{Proposition}\label{prop:rec_dualpolar}
(i) If $0 \le i \le d$, then
$C_j^{d+1}(i) = q^{d+e-i} C_{j-1}^{d}(i) + C_j^{d}(i)$.

(ii) If $1 \le i \le d+1$, then
$C_j^{d+1}(i) = -q^{i-1} C_{j-1}^{d}(i-1) + C_j^{d}(i-1)$. \qed
\end{Proposition}

\noindent
Since these two values are equal, one also has
$C_j(i-1) = C_j(i) + q^{i-1} C_{j-1}(i-1) + q^{d+e-i} C_{j-1}(i)$.

\medskip
%
%
We have
$C_1(i) = q^e \qnom{d-i}{1} - \qnom{i}{1}$ and
$C_d(i) = (-1)^i q^{\binom{d}{2}+(d-i)(e-i)}$,
and see that for $j = 1$ and for $j = d$ the sequence
$|C_j(i)|$ $(0 \le i \le d)$ is unimodal, with smallest element
$|C_j(i)|$ for $i = \lfloor(d+e+1)/2\rfloor$, largest element
$C_j(0)$ and second largest element $|C_j(d)|$ if $e \le 1$,
and $|C_j(1)|$ if $e > 1$.
This is what we try to prove for all $j$.

There are small exceptions. E.g. for $(q,d,e) = (2,5,1)$ the
$j = 4$ column of $P$ is not unimodal, and the $j = 2$ column
has its minimum earlier:
$$
P = \left(\begin{array}{cccccc}
1 & 62 & 1240 & 9920 & 31744 & 32768\\
1 & 29 & 250 & 680 & 64 & -1024\\
1 & 11 & 16 & -76 & -80 & 128\\
1 & -1 & -20 & 20 & 64 & -64\\
1 & -13 & 40 & 20 & -176 & 128\\
1 & -31 & 310 & -1240 & 1984 & -1024\\
\end{array}\right) .
$$
More generally, if $(q,e) = (2,1)$, then
$|C_{d-1}(2)| > |C_{d-1}(1)| = q^{\binom{d-1}{2}}$ for all $d \ge 2$,
and the sequence $|C_{d-1}(i)|$ is not unimodal for $(q,e) = (2,1)$, $d \ge 5$.
%
%

For $e=1$ we have the coincidence $|C_d(1)| = |C_d(d)|$.
More generally, $|C_d(i)| = |C_d(d+e-i)|$ for integral $e$
and $e \le i \le d$.

For $e=0$ the graphs $C_q(d,e)$ are bipartite, and we have
$C_j(d-i) = (-1)^j C_j(i)$.

\subsection{The smallest eigenvalue}

The following conjecture is a variation of Lemma 47 in \cite{FIKM}
where the authors investigated the sum of the relations
$\{ d-j, d-j+1, \ldots, d \}$ instead of just the $j$th relation.

\begin{Conjecture}
The sequence $|C_j(i)|$ ($j$ fixed, $0 \le i \le d$) is unimodal
if not $(q,e) = (2,1)$ and not $(q,e,j) = (2,2,d-4)$, $8 \le d \le 12$.
If it is unimodal with minimum at $i_0$, and $i_1 = \lfloor(d+e+1)/2\rfloor$,
then $i_0 = i_1$ for $e = 0,\frac12,\frac32$, and
$|i_0 - i_1| \le 1$ for $e=1,2$,
except that $i_0 = i_1-2$ for $(q,e,j,d) = (2,1,3,4),(2,2,3,7)$.
\end{Conjecture}
%
%

\begin{Conjecture}
The index $i_{\rm min}$ of the smallest among the $C_j(i)$
($j$ fixed, $0 \le i \le d$) is
$$
i_{\rm min} = \left\{
\begin{array}{ll}
1 & \mbox{if $j = d$ and ($j$ is even or $e \ge 1$)} \\
d & \mbox{if $j$ is odd and ($j < d$ or $e \le 1$)} \\
\lfloor (d-j+2)/2 \rfloor & \mbox{if $j$ is even, $e = 0$} \\
(d-j+2)/2 & \mbox{if $j$ and $d$ are even, $e=\frac12$ or $e=1$} \\
(d+j-1)/2 & \mbox{if $j$ is even, $d$ is odd, $e=\frac12$ or $e=1$} \\
(d+j)/2 & \mbox{if $j$ and $d$ are even, $e=\frac32$ or $e=2$} \\
(d-j+3)/2 & \mbox{if $j$ is even, $d$ is odd, $e=\frac32$ or $e=2$} \\
\end{array}\right.
$$
except that when $q=2$ and $e=2$ and $d$ is even and $j \ge d-4$
one finds $i_{\rm min} = 2$ for $j = d-2$, $d \ge 6$
and $i_{\rm min} = 3$ for $j = d-4$, $d \ge 14$.
\end{Conjecture}

We show the second case of this conjecture in
Corollary \ref{dualpolar-smallest}.
We can show the conjecture for some more cases if $q \geq 11$,
but omit the details.

\begin{Proposition} Let $1 \leq j \leq d$.

(i) $C_j(1) < 0$ if and only if $j=d$ or $(j, e) = (d-1, 0)$.

(ii) Let $d \geq 3$. Then $|C_j(2)| \leq |C_j(1)|$
    unless $(q,j,e) = (2,d-1,1)$.

(iii) Let $1 \leq i \leq d$. Then $|C_j(i)| \leq |C_j(d)|$ if $i \geq 2$
    or $e \leq 1$.

(iv) $|C_j(1)| \leq |C_j(d)|$ if $e \leq 1$ with equality only if $(j,e)=(d,1)$.
\end{Proposition}
\Proof
  (i) This is immediate from $C_j(1) = q^{\binom{j}{2}+je} \qnom{d-1}{j}
       - q^{\binom{j-1}{2}+(j-1)e} \qnom{d-1}{j-1}$.
  
  (ii) We can assume $1 < j < d$ as we did already show the claim for $j=1$
      and $j=d$. Rename $d$ to $d+1$, so that $d \ge 2$ and $j \le d$.
      We have
      $C^{d+1}_{j}(1) = q^{\binom{j}{2}+je} \qnom{d}{j} - q^{\binom{j-1}{2}+(j-1)e} \qnom{d}{j-1}$,
      and $C^{d+1}_{j}(2) = C^{d}_j(1) - q C^{d}_{j-1}(1)$
      by Proposition \ref{prop:rec_dualpolar} (ii).
      Dividing the expression $|C^{d+1}_{j}(2)| \leq |C^{d+1}_{j}(1)|$ by $q^{\binom{j-1}{2}+(j-1)e} \qnom{d-1}{j-1}$ and 
      simplifying yields the claim.

  (iii) Note that $C_j(d) = (-1)^{j} q^{\binom{j}{2}} \qnom{d}{d-j}$
        has alternating sign.
	Use induction on $d$.
        By Proposition \ref{prop:rec_dualpolar} (i) and (ii), 
	\begin{align*}
	  |C_j^{d+1}(i)| &= |q^{d+e-i} C^{d}_{j-1}(i) + C^{d}_j(i)| \\
	  &\leq |q^d C^{d}_{j-1}(d)| + |C^{d}_{j}(d)| = |C^{d+1}_j(d+1)|.
	\end{align*}

  (iv) This is immediate from the expressions for $C_j(1)$ and $C_j(d)$.
\qed

\begin{Corollary} \label{dualpolar-smallest}
%
%
  Let $d \geq 3$ and $1 \leq j \leq d$.
  Then
  
  (i) $|C_j(1)| = \max\{ |C_j(i)|: 1 \leq i \leq d )$ if $e > 1$ or $(j,e) = (d,1)$.
  
  (ii) $|C_j(d)| = \max\{ |C_j(i)|: 1 \leq i \leq d )$ if $e \leq 1$.
  
  (iii) If $j < d$ is odd, then $C_j(d) = \min\{ C_j(i): 0 \leq i \leq d \}$.
\end{Corollary}
\Proof
  We only have to show (iii). Here we only have to show that $C_j(d)$ is negative.
  This follows from Proposition \ref{lastrow}.
\qed

\section{Bilinear forms graphs}
The bilinear forms graphs $H_q(d,e)$ are the graphs with as vertices
$d \times e$ matrices over $\Ff_q$, adjacent when the difference
has rank 1. W.l.o.g. we assume $d \le e$.
The eigenmatrix $P$ has entries $P_{ij} = B_j(i)$, where
$$
B_j(i) = \sum_{h=0}^j (-1)^{j-h} \, q^{eh+\binom{j-h}{2}} \,
\qnom{d-h}{d-j} \qnom{d-i}{h}
$$
(Delsarte \cite{Delsarte78}, Theorem A2).

The valencies here are
$k_j = B_j(0) = \bnom{d}{j}\bnom{e}{j} \prod_{h=1}^{j} (q^j - q^{j-h})$
(\cite{BCN}, p. 281).

The eigenvalues of $H_q(d,e)$ are $\theta_i = (q^{d+e-i}-q^d-q^e+1)/(q-1)$.

The scheme is self-dual, so that $P_{ij}/P_{0j} = P_{ji}/P_{0i}$,
and $P_{ij}$ and $P_{ji}$ have the same sign.

\subsection{Identities}
Let us write $B_j^{d,e}(i)$ instead of $B_j(i)$ when it is necessary
to make the dependence on $d$ and $e$ explicit.

\begin{Proposition} {\rm (Delsarte \cite[Proof of Theorem A2]{Delsarte78})}
$$
B_j^{d,e}(i) - B_j^{d,e}(i+1) = q^{d+e-i-1} B_{j-1}^{d-1,e-1}(i) . 
$$
\end{Proposition}

\begin{Proposition} {\rm (Stanton, \cite[Prop. 1(ii),(iii)]{Stanton81})}

(i)
$(q^{d-j+1}-1) B_j^{d+1,e}(i) = (q^{d+1}-q^i) B^{d,e}_j(i) + (q^i-1)B^{d,e}_j(i-1)$.

(ii)
$(q^{e-j+1}-1) B_j^{d,e+1}(i) = (q^{e+1}-q^i) B^{d,e}_j(i) + (q^i-1)B^{d,e}_j(i-1)$.
\end{Proposition}


\subsection{Negative \texorpdfstring{$B_j(1)$}{Bj(1)}}
For the bilinear forms graphs the $i=1$ row of $P$ has only a single
negative value.

\begin{Proposition}\quad \label{bilin-neg}

(i) $B_j(1) < 0$ if and only if $j = d$, and otherwise $B_j(1) > 0$.

(ii) $B_d(1)$ is the smallest eigenvalue of the distance-$d$ graph,
and the second largest in absolute value.
\end{Proposition}

\Proof
(i) This follows from
$B_1(i) = (q^{d+e-i} - q^d - q^e + 1)/(q-1)$
and $e \ge d$ and the fact that $B_1(i)$ and $B_i(1)$ have the same sign.

(ii) Proposition \ref{lastrow} gives
$B_j(d) = (-1)^j \bnom{d}{j} q^{\binom{j}{2}}$, and it follows that
$B_d(i) = (k_d/k_i) (-1)^i \bnom{d}{i} q^{\binom{i}{2}}$.
The claim follows using
$k_i = \bnom{d}{i}\bnom{e}{i} q^{\binom{i}{2}} \prod_{h=1}^{i} (q^h - 1)$.
\qed

\begin{Lemma} \label{lemma:bilinbds}
Let $1 \le j \le d-1$ and either $j \le d-2$ or $q > 2$ or $q=2, e>d$.
Then $|B_j(2)| \le |B_j(1)|$. If $j=d-1$ and $q=2$ and $e=d$, then
$|B_j(2)|/|B_j(1)| = (2^{d-1}+1)/(2^{d-1}-1)$.
\end{Lemma}
\Proof
Find $B_j(1)$ and $B_j(2)$ from $B_1(i)$ and $B_2(i)$ and the relation
$P_{ij}/k_j = P_{ji}/k_i$. Abbreviate $q^n-1$ with $[n]$.
One gets
$$
\frac{B_j(2)}{B_j(1)} =
\frac{ q[d][d-1] - (q+1)q^e[d-1][d-j] + q^{2e}[d-j][d-j-1] }
     { q[d-1][e-1](q^{d+e-j}-q^d-q^e+1) } .
$$
The numerator is of the form $A-B+C$ where $B \ge A \ge 0$ and $C \ge 0$.
If $j \le d-3$, or $j=d-2$, $q>2$, or $j=d-2$, $q=2$, $e>d$, then $C \ge B$.
Now estimate the numerator with $C$ and find that
$|\frac{B_j(2)}{B_j(1)}| \le 1$. The same conclusion follows
by direct computation in the case $j=d-2$, $q=2$, $e=d$.
This leaves the case $j = d-1$ (with $C = 0$).
Again treat the cases $q > 2$ and $q=2$, $e > d$ separately
and find the same conclusion.
\qed

\bigskip
As the scheme is self-dual, so that $P_{ij}/k_j = P_{ji}/k_i$,
the recurrence
$c_{j+1} P_{i,j+1} \maysplit = (\theta_i - a_j) P_{ij} - b_{j-1} P_{i,j-1}$
implies 
$b_i P_{i+1,j} = (\theta_j - a_i) P_{ij} - c_i P_{i-1,j}$.
In our case this gives (after multiplication by $q-1$)
\begin{align*}
& q^{2i} [d-i] [e-i] B_j(i+1) \\
&  = \left( q^e [d-j] - [d] - [i] (q^e+q^d-q^i-q^{i-1}-1) \right) B_j(i)
   - q^{i-1} [i] B_j(i-1) ,
\end{align*}
again with the abbreviation $[n] = q^n-1$.

\begin{Theorem} \label{bilinear-2ndlargest}
For $q \ge 4$, $|B_j(1)| \geq |B_j(i)|$ for
$1 \le i \le d$, $0 \le j \le d$.
\end{Theorem}
\Proof
For $j=0$ the claim is trivial, so we assume $j \ge 1$.
By Propostion \ref{bilin-neg} we can assume $j < d$.
Now $|B_j(i)| \le |B_j(1)|$ follows by induction on $i$,
starting with Lemma \ref{lemma:bilinbds} for $i = 2$, and
using the recurrence for $i > 2$.
We have to show that
$\max(|q^e [d-j] - [d]|,~[i] (q^e+q^d-q^i-q^{i-1}-1)) +
q^{i-1} [i] \le q^{2i} [d-i] [e-i]$,
and that is easily checked, assuming $q \ge 4$.
\qed

\begin{Conjecture}
For $q \ge 3$, or $q = 2$ and $d \neq e$,
$B_j(d-j+1)$ is the smallest eigenvalue in the distance-$j$ graph
for $1 \le j \le d$.
\end{Conjecture}

Let $b_{i,j}(h)$ be the exponent of $q$ in the $h$-th term of the
expression for $B_j(i)$ if we approximate $\qnom{n}{k}$ with $q^{k(n-k)}$.
That is, let
\begin{align*}
b_{i,j} (h) = h(d+e-i-h) + (d-j)(j-h) + \binom{j-h}{2} .
\end{align*}
Let $h_0 = e-i+\frac12$.
Then the quadratic expression $b_{i,j}(h)$ is maximal for $h = h_0$,
and $b_{i,j}(h_0+x) = b_{i,j}(h_0) - \frac12 x^2$.
Let $h_{\rm max} = \min(j,d-i)$.
The terms occurring in the sum have indices $h$ with $h \le h_{\rm max} < h_0$,
so the term with largest index has largest exponent.

\begin{Lemma}\label{lem:bad_bnd_Bij}
  Let $q \geq 4$ and put $s := b_{i,j}(h_{\rm max})$. We have
  \begin{align*}
    \frac{5}{9} q^s < |B_j(i)| < \frac{13}{4} q^s.
  \end{align*}
\end{Lemma}
\Proof
The  expression for $B_j(i)$ is an alternating series
with terms decreasing in absolute value after the first,
so we can estimate $B_j(i)$ by the main term with an error
not larger than the second term.
\qed

\begin{Proposition} Let $q \geq 4$.
The sign of $B_j(i)$ is $(-1)^{\max(0,i+j-d)}$.
The smallest among the $B_j(i)$ for fixed $j$ is $B_j(d-j+1)$.
\end{Proposition}
\Proof
The sign of $B_j(i)$ is that of the main term.
The negative terms are $B_j(d-j+1+2t)$.
Increasing $i$ by 2 (from $d-j+1+2t$ to $d-j+3+2t$)
means decreasing $s$ by at least $2$, and
since $\frac{5}{9}q^2 > \frac{13}{4}$ that decreases
the absolute value. So $B_j(d-j+1)$ is most negative.
\qed

\section{Alternating forms graphs}
The alternating forms graphs $A_q(n)$ are the graphs with as vertices
the skew symmetric $n\times n$ matrices over $\Ff_q$ with zero diagonal,
adjacent when the difference has rank 2.

Let $d = \lfloor n/2 \rfloor$.
The graph $A_q(n)$ is distance-regular with diameter $d$.
The eigenmatrix $P$ has entries $P_{ij} = A_j(i)$, where
$$
A_j(i) = \sum_{h=0}^j (-1)^{j-h} \, q^{(j-h)(j-h-1)} \, q^{hm} \,
\bnom{d-h}{d-j}_b \bnom{d-i}{h}_b.
$$
Here the Gaussian coefficients have base $b = q^2$
and $m = n(n-1)/(2d) = 2n-2d-1$ so that $\{m,2d\} = \{n-1,n\}$
and $m$ is odd (Delsarte \cite[(15)]{Delsarte75}).

The valencies here are $k_j = A_j(0) = q^{j(j-1)} \,
\prod_{i=0}^{2j-1} (q^{n-i}-1) / \prod_{i=1}^j (q^{2i}-1)$.

The eigenvalues of $A_q(n)$ are
$\theta_i = (q^{2n-2i-1}-q^n-q^{n-1}+1)/(q^2-1)$.

The scheme is self-dual, so that $P_{ij}/P_{0j} = P_{ji}/P_{0i}$,
and $P_{ij}$ and $P_{ji}$ have the same sign.

\subsection{Identities}

Let us write $A_j^{n}(i)$ instead of $A_j(i)$ when it is necessary
to make the dependence on $n$ explicit.

\begin{Proposition} \label{alt-recur}
{\rm (Delsarte \cite[(66)]{Delsarte75})}
$A_j^{n}(i) = A_j^{n}(i{-}1) - q^{2n-2i-1} A_{j-1}^{n-2}(i{-}1)$.
\end{Proposition}

\begin{Proposition} \label{alt-lastcol}
$A_d(i) = - (q^{m-2i} - 1) A_d(i+1)$ for $0 \le i \le d-1$. \qed
\end{Proposition}

\subsection{The smallest and the second largest eigenvalue}

We determine the smallest eigenvalue, and the second largest in absolute value,
for the distance-$j$ graphs of $A_q(n)$.

\begin{Theorem} \label{alt-2ndlargest} Let $1 \le j \le d$.

(i) $\min_{0\leq i\leq d} A_j(i) = A_j(d-j+1)$.

(ii) $\max_{1\leq i \leq d} |A_j(i)| = |A_j(1)|$.

(iii) Let $0 \le i \le d-1$, $1 \le j \le d$. Then:

a) $|A_j(i)| < |A_j(i+1)|$ if and only if $(q,n,i) = (2,2d,d-1)$
and $1 \le j \le d-1$.

b) $|A_j(i)| = |A_j(i+1)|$ if and only if $(q,n,i) = (2,2d,d-1)$
and $j = d$.

c) In all other cases $|A_j(i)| > |A_j(i+1)|$.
\end{Theorem}

The proof of this theorem is given below.

For $(q,n) = (2,4)$ we have
$$
P = \left(\begin{array}{ccc}
1 & 35 & 28\\
1 & 3 & -4\\
1 & -5 & 4\\
\end{array}\right) .
$$

\medskip
Let $a_{i,j}(h)$ be the exponent of $q$ in the $h$-th term of the
expression for $A_j(i)$ if we approximate $\bnom{n}{k}_b$ with $q^{2k(n-k)}$.
Then
\begin{align*}
a_{i,j}(h) &= (j-h)(j-h-1) + hm + 2(d-j)(j-h) + 2h(d-i-h) \\
 &= -h^2 + h(m+1-2i) + j(2d-j-1) .
\end{align*}
This quadratic function of $h$ is maximal for $h_0 = \frac{m+1}{2} - i$.
The nonzero terms in the expression for $A_j(i)$ have indices $h$
with $0 \le h \le \min(d-i,j)$. Since $h_0 = d-i$ if $n$ is even,
and $h_0 = d-i+1$ if $n$ is odd, the term with the largest exponent
is the one with index $\min(d-i,j)$.

\begin{Proposition}
If $i+j \le d$, then
$$
0 \le 1 - \frac{A_j(i)}{q^{jm} \bnom{d-i}{j}_b} <
\frac{2}{q^{m+2-2i-2j}} .
$$
In particular, $A_j(i) > 0$.
\end{Proposition}
\Proof
Use that $a_{i,j}(h_0-x) = a_{i,j}(h_0) - x^2$.
If $i+j \le d$, then $\min(d-i,j) = j$.
The sum is alternating, and since $b = q^2 \ge 4$
and $(1+2q^{-2})^2 < q^3$ it follows from Lemma \ref{gauss_bnds} (iii,v)
that terms after the first (reading down from largest $h$)
decrease in size, and the difference between $A_j(i)$
and the first term is not larger than the second term.
(That is, $A_j(i) = T_0 - T_1 + T_2 - \cdots$ where all $T_{\ell}$ have the same
sign, and $|T_1| \ge |T_2| \ge \cdots$.
Our conclusion will be $A_j(i) = T_0 - \gamma T_1$ with $0 \le \gamma \le 1$,
that is, $1 - \frac{T_1}{T_0} \le \frac{A_j(i)}{T_0} \le 1$.)
Estimate the absolute value of second term divided by the first,
using Lemma \ref{gauss_bnds} (ii), by
$$
\frac{\bnom{d-j+1}{1} \bnom{d-i}{j-1}}{q^m \bnom{d-i}{j}} =
q^{-m} \frac{b^{d-j+1}-1}{b-1} \frac{b^j-1}{b^{d-i-j+1}-1} <
q^{-m+2i+2j-2}\frac{q^4}{(q^2-1)^2}.
$$
If $n$ is odd, $m+2-2i-2j = 2d+3-2i-2j \ge 3$.
If $n$ is even, $m+2-2i-2j = 2d+1-2i-2j \ge 1$.
In both cases, the RHS of the inequality is less than 1.
\qed

\begin{Proposition} \label{alt-below}
If $s := i+j-d \ge 0$, then
$$
0 \le 1 - \frac{A_j(i)}{(-1)^s q^{s(s-1)+(d-i)m} \bnom{i}{d-j}}
 \le \frac{\bnom{i+1}{d-j} \bnom{d-i}{1}}{q^{m-2s} \bnom{i}{d-j}}
 < \frac{q^3}{(q^2-1)^2} \frac{1}{q^{2n-4d}} < 1 .
$$
In particular, $A_j(i)$ has sign $(-1)^s$.
\end{Proposition}
\Proof
If $i+j \ge d$, then $\min(d-i,j) = d-i$.
Again the difference between $A_j(i)$ and the first term is not larger
than the second term.
Estimate the absolute value of second term divided by the first,
using Lemma \ref{gauss_bnds} (ii), by
$$
\frac{\bnom{i+1}{d-j} \bnom{d-i}{1}}{q^{m-2s} \bnom{i}{d-j}} =
\frac{(b^{d-i}-1)(b^{i+1}-1)}{q^{m-2s}(b-1)(b^{s+1}-1)} <
\frac{q^{-m+2d+2}}{(q^2-1)^2}.
$$
Finally, $m-2d = 2n-4d-1$.
\qed

\medskip
\Proofn {\bf of Theorem \ref{alt-2ndlargest}.}
First of all, by Proposition \ref{alt-lastcol} all statements are true
for $j = d$, so we may suppose $1 \le j \le d-1$.

Next, prove part (iiic). We have 
$A^n_j(i+1) = A^n_j(i) - q^{2n-2i-3} A^{n-2}_{j-1}(i)$.

If $i+j+1 \le d$, then each of
$A^n_j(i+1)$, $A^n_j(i)$, $A^{n-2}_{j-1}(i)$ is positive,
and $|A^n_j(i+1)| < |A^n_j(i)|$ follows from
$0 < A^n_j(i+1) < A^n_j(i)$.


If $i+j \ge d$, use the (strong form of the) second proposition to find
\begin{align*}
\left|\frac{A_j(i+1)}{A_j(i)}\right| &<
\frac{q^{2s-m}\,\qnom{i+1}{d-j}}
     {\left(1-\txtfrac{\bnom{i+1}{d-j} \bnom{d-i}{1}}
                 {q^{m-2s} \bnom{i}{d-j}}\right)\,\qnom{i}{d-j}} =
\frac{1}{q^{m-2s} \, \txtfrac{b^{s+1}-1}{b^{i+1}-1} -
                  \txtfrac{b^{d-i}-1}{b-1}} \\
&= \frac{(b-1)(b^{i+1}-1)}
        {q^{m-2s}(b-1)(b^{s+1}-1)-(b^{d-i}-1)(b^{i+1}-1)} \\
&< \frac{(b-1)b^{i+1}}{(b-1)^2 q^m-b^{d+1}}.
\end{align*}
If $n$ is odd, then $m = 2d+1$, and the RHS is at most $\frac{3}{14}$
(since $i \le d-1$ and $q \ge 2$). If $n$ is even, then $m = 2d-1$.
Now if $q \ge 3$ then the RHS is at most $\frac{24}{37}$.
If $q = 2$ and $i \le d-3$ then the RHS is at most $\frac{3}{8}$.
For $q = 2$ and $i = d-2$ we use the sharper form of the last inequality.
The claim $|A_j(i+1)/A_j(i)| < 1$ follows from
$q^{m-2s} \, \txtfrac{b^{s+1}-1}{b^{i+1}-1} - \txtfrac{b^{d-i}-1}{b-1} > 1$,
which is true since $b^d > b^{d-s} - 12$.
That proves part (iiic).

Part (iiib) is the case $(q,n,i) = (2,2d,d-1)$
of Proposition \ref{alt-lastcol}.

Part (iiia) follows from
$-\frac{A_j(d-1)}{A_j(d)} = q^{m-2j+2} \frac{b^j-1}{b^d-1} - 1$.
This is larger than 1, unless $q = 2$ and $n$ is even.

That proves part (iii). Now part (ii) follows, except in the case
$(q,n) = (2,2d)$. We show that in this case
$|A_j(d-2)| > |A_j(d)|$. Indeed,
$|A_j(d)| = q^{j(j-1)} \bnom{d}{d-j}_b$ and
$|A_j(d-2)| > (1-\gamma) q^{(j-2)(j-3)+2m}\bnom{d-2}{d-j}_b$, where
$\gamma < \frac{8}{9}$ and the desired inequality follows from
Lemma \ref{gauss_bnds} (iii),(v).

Finally part (i) follows, since the smallest among the $A_j(i)$
is the first one that is negative.
\qed

\section{Hermitian forms graphs}

The Hermitian forms graphs $Q_q(d)$ are the graphs with as vertices
the Hermitian $d\times d$ matrices over $\Ff_{q^2}$,
adjacent when the difference has rank 1.

The graph $Q_q(d)$ is distance-regular with diameter $d$. The eigenmatrix $P$
has entries $P_{ij} = Q_j(i)$, where
$$
Q_j(i) = (-1)^j \sum_{h=0}^j \, (-q)^{\binom{j-h}{2} + hd} \, \bnom{d-h}{d-j}_b \bnom{d-i}{h}_b.
$$
Here the Gaussian coefficients have base $b=-q$.
This formula was taken from Schmidt \cite{Schmidt2017}.
An expression in terms of $q$-Krawtchouk polynomials was
given in Stanton \cite{Stanton81}.

The eigenvalues of $Q_q(d)$ are
$\theta_i = ((-q)^{2d-i} - 1)/(q+1)$.

The scheme is self-dual, so that $P_{ij}/P_{0j} = P_{ji}/P_{0i}$,
and $P_{ij}$ and $P_{ji}$ have the same sign.

\subsection{Identities}

Let us write $Q_j^d(i)$ instead of $Q_j(i)$ when it is necessary
to make the dependency on $d$ explicit.

\begin{Proposition} {\rm (\cite[Lemma~7]{Schmidt2017})} \label{prop:rec_hermitian}
 $Q_j^{d}(i) = Q_j^d(i-1) + (-q)^{2d-i} Q_{j-1}^{d-1}(i-1)$.
\end{Proposition}


\subsection{The smallest and the second largest eigenvalue}

\begin{Conjecture}
  (i) If $j$ is odd, then $Q_j(1) \leq Q_j(i)$ for $0 \leq i \leq d$.
  
  (ii) If $j$ is even, $j \ge 2$, then $Q_j(d-j+2) \leq Q_j(i)$ for $0 \leq i \leq d$.
\end{Conjecture}

\begin{Conjecture}
  Let $d \geq 3$. Then $|Q_j(i)| < |Q_j(1)|$ for $2 \leq i \leq d$.
\end{Conjecture}
In the following we prove both conjectures for $q \geq 4$.

Let $q_{i,j}(h)$ be the exponent of $q$ in the $h$-th term of the expression for $Q_j(i)$
if we approximate $|\bnom{n}{k}_{-q}|$ with $q^{k(n-k)}$. Then
\begin{align*}
  q_{i,j}(h) &= (d-j)(j-h) + h(d-i-h) + (j-h)(j-h-1)/2 + hd.
\end{align*}
Let $h_0 = d - i + \frac{1}{2}$. Then the quadratic expression $q_{i,j}(h)$ is maximal
for $h = h_0$, and $q_{i,j}(h_0+x) = q_{i,j}(h_0) - \frac{1}{2}x^2$. 
Let $h_{\max} = \min(j, d-i)$. The terms occurring in the sum have indices $h$ 
with $h \leq h_{\max} < h_0$, so the term with the largest index has the largest exponent.

\begin{Proposition}\label{prop:estimate_hermitian_form} Let $d \ge 2$,
  $j \ge 1$ and $q \geq 4$.
  Set $S = S(i) := \bnom{d-i}{j} (-q)^{jd}$ if $d-i \geq j$
  and $S = S(i) := \bnom{i}{d-j} (-q)^{\binom{i+j-d}{2} + (d-i)d}$ otherwise. Then
  $$
    \left| Q_j(i) - (-1)^j S \right| \le \frac{11}{27} |S|.
  $$
  In particular, the sign of $Q_j(i)$ is the sign of $(-1)^j S$.
\end{Proposition}
\Proof
  If we divide the absolute value of the $h$-th term in the expression
  by the absolute value of the $(h-1)$-th term in the expression,
  then we obtain, using $1 \le h \leq \min(j, d-i)$, and (for $m > 0$)
$$\begin{array}{cccl}
  (1-q^{-m})q^m & \le |b^m-1| \le & q^m & \mbox{~~if $m$ is even,} \\
  q^m & \le |b^m-1| \le & (1+q^{-m})q^m & \mbox{~~if $m$ is odd,}
\end{array}$$
  and $q \ge 4$, that
  \begin{align*}
    \left| \frac{b^{j-h+1}-1}{b^{d-h+1}-1}  \cdot \frac{b^{d-i-h+1}-1}{b^{h}-1}
           \cdot b^{h-j+d} \right| 
    \geq \frac{(1-q^{-2})^2}{(1+q^{-1})(1+q^{-3})} \, q^{d-i-h+1}
    \geq \frac{9}{13} q^a
  \end{align*}
  where $a = 1$ if $h = d-i$, and $a = 2$ otherwise.
  Then (again using $q \ge 4$)
  \begin{align*}
     \left| Q_j(i) - (-1)^j S \right| \leq
      \frac{13}{36} \sum_{h\ge0} 10^{-h} \, |S| \leq \frac{11}{27} |S|.
  \end{align*}
  This shows the assertion.
\qed

\begin{Theorem} Let $j \geq 1$ and $q \geq 4$. \label{herm-smallest-2ndlargest}

  (i) Let $d \geq 3$. Then $|Q_j(i+1)| < |Q_j(i)|$ for $0 \leq i \leq d-1$.
  
  (ii) If $j$ is odd, then $Q_j(1) \leq Q_j(i)$ for $0 \leq i \leq d$.
  
  (iii) If $j$ is even, then $Q_j(d-j+2) \leq Q_j(i)$ for $0 \leq i \leq d$.
\end{Theorem}
\Proof 
  (i) 
  By Proposition \ref{prop:estimate_hermitian_form}, we have
  $|Q_j(i)| \geq \frac{16}{27} |S(i)|$ and $|Q_j(i+1)| \leq \frac{38}{27} |S(i+1)|$.
  We have to show that $|S(i)|/|S(i+1)| > \frac{19}{8}$.
  If $i+j \le d-1$,
  $$
    \frac{|S(i)|}{|S(i+1)|}
     = \left| \frac{\bnom{d-i}{j}}{\bnom{d-i-1}{j}} \right|
     = \left| \frac{b^{d-i}-1}{b^{d-i-j}-1} \right|
     > \frac{1-q^{-1}}{1+q^{-1}} \, q^j > \frac{19}{8}.
  $$
  If $i+j \ge d$,
  $$
    \frac{|S(i)|}{|S(i+1)|}
     = \left| \frac{\bnom{i}{d-j}}{\bnom{i+1}{d-j}} b^{-i-j+2d} \right|
     = \left| \frac{b^{i+j-d+1}-1}{b^{i+1}-1} \, b^{-i-j+2d} \right|
     > \frac{1-q^{-1}}{1+q^{-1}} \, q^{d-i} > \frac{19}{8}.
  $$
  
  (ii) and (iii) 
  By Proposition \ref{prop:estimate_hermitian_form} and part (i),
  we only have to find the smallest $i$ for which $(-1)^j S(i)$ is negative.
  The sign of $(-1)^j \bnom{d-i}{j} b^{jd}$ is positive for $j$ even,
  and $(-1)^{jd + d-i} = (-1)^i$ for $j$ odd. This proves part (ii).
  The sign of $(-1)^j \bnom{i}{d-j} (-q)^{\binom{i+j-d}{2} + (d-i)d}$
  where $j$ is even, is $(-1)^{\binom{i+j-d}{2}}$,
  hence is positive for $i=d-j+1$ and negative for $i=d-j+2$.
  This shows (iii).
\qed

\section*{Acknowledgements}
We thank Sasha Barg for bringing reference \cite{DK} to our attention.
The research of Sebastian M. Cioab\u{a} and Matt McGinnis is supported by NSF grant DMS-1600768.
The research of Ferdinand Ihringer is supported by ERC advanced grant 320924 and a postdoctoral fellowship of the Research Foundation - Flanders (FWO).

\end{document}